\documentclass[format=acmsmall,manuscript,screen]{acmart}
\pdfoutput=1

\newif\ifarxiv
\arxivfalse

\setcopyright{none}

\usepackage[utf8]{inputenc}
\usepackage{array}
\usepackage{algorithm}
\usepackage{algpseudocode}
\usepackage{amsmath}
\usepackage{booktabs}
\usepackage[english]{babel}
\usepackage{enumitem}
\usepackage{hyperref}
\usepackage{pgfplots}
\usepackage{pgfplotstable}
\usepackage{cleveref}
\usepackage{tikz}
\usetikzlibrary{calc}
\usepackage{verbatim}
\usepackage{subcaption}
\usepackage{xspace}

\usetikzlibrary{shapes,arrows,positioning,decorations.pathreplacing,
        intersections,matrix,patterns,fadings,calc,decorations.markings,arrows.meta}

\pgfplotsset{compat=1.17}

\newcommand{\bfZ}{\mathbf{Z}}
\newcommand{\bfY}{\mathbf{Y}}
\newcommand{\bfy}{\mathbf{y}}
\newcommand{\bfk}{\mathbf{k}}
\newcommand{\bfkappa}{\boldsymbol{\kappa}}
\newcommand{\bff}{\mathbf{f}}
\newcommand{\bfF}{\mathbf{F}}
\newcommand{\bfc}{\mathbf{c}}
\newcommand{\bfb}{\mathbf{b}}

\newcommand{\Irksome}{\texttt{Irksome}}
\newcommand{\RK}{Runge--Kutta}
\newcommand{\RKN}{Runge--Kutta--Nystr\"om}

\usepackage{listings}
\definecolor{DarkBlue}{rgb}{0.00,0.00,0.55}
\definecolor{DarkRed}{rgb}{0.55,0.00,0.00}
\definecolor{DarkGreen}{rgb}{0.00,0.55,0.00}
\definecolor{Bittersweet}{rgb}{1.0, 0.44, 0.37}
\definecolor{Purple}{rgb}{0.5, 0.0, 0.5}

\acmJournal{TOMS}

\title{Automated \RKN{} time stepping for finite element methods in \Irksome}

\author{Robert C. Kirby}
\affiliation{%
  \institution{Baylor University}
  \department{Department of Mathematics}
  \streetaddress{1410 S.~4th St.}
  \city{Waco}
  \state{TX}
  \country{USA}}
\email{robert_kirby@baylor.edu}

\author{Scott P. MacLachlan}
\affiliation{%
  \institution{Memorial University of Newfoundland}
  \department{Department of Mathematics and Statistics}
  \city{St.~John's}
  \state{NL}
  \country{Canada}
}
\email{smaclachlan@mun.ca}

\author{Pablo D. Brubeck}
\affiliation{%
  \institution{University of Oxford}
  \department{Mathematical Institute}
  \city{Oxford}
  \country{UK}}
\email{brubeckmarti@maths.ox.ac.uk}

\numberwithin{equation}{section}
\crefname{algorithm}{Algorithm}{Algorithms}
\crefname{figure}{Fig.}{Figs.}
\crefname{table}{Table}{Tables}

\begin{abstract}
  \texttt{Irksome} is a library based on the Unified Form Language (UFL) that
  automates the application of Runge--Kutta time-stepping methods for finite element spatial discretizations of partial differential equations (PDEs).
  This paper describes recent updates to \texttt{Irksome} that allow users to express semidiscrete forms of PDEs that contain second-order temporal derivatives,
  whence it generates stage-coupled variational problems to be solved at each time step for \RKN{} methods.
  Firedrake then generates code for these variational problems and provides a rich interface to PETSc for solving them.
  Directly discretizing second-order time derivatives with \RKN{} methods provides several advantages relative to discretizing a rewritten first-order system with a standard \RK{} method.
  Besides working with an interface closer to the problem formulation in UFL, avoiding these auxiliary variables means that \RKN{} methods lead to smaller algebraic systems and better run-time.
  Our numerical results indicate that, with effective preconditioning, fully implicit \RKN{} methods can be made competitive with more traditional explicit methods for wave equations.  They are also (essentially) required to discretize wave-type equations with higher-order spatial derivatives.  We also provide numerical experiments for fully dynamic poroelasticity, a system of mixed temporal order, where our time-stepping and algebraic solvers perform effectively even as we approach the incompressible limit.  
\end{abstract}


\begin{document}
\maketitle

\section{Introduction}

The Unified Form Language, implemented by tools such as Firedrake~\cite{FiredrakeUserManual, Rathgeber:2016}, FEniCS~\cite{Logg:2012}, and Dune~\cite{dedner2020python} provides a high-level Python library for describing finite element discretizations and solutions of partial differential equations (PDEs).
For example, a Firedrake implementation of the standard Poisson equation is ably described in Figure~\ref{fig:fdpoisson}.
Within the call to \lstinline{solve}, code is generated and compiled to assemble the relevant stiffness matrix and load vector, and the solution process is performed in PETSc (steered beyond some default settings via an optional keyword argument).

\begin{figure}
  \begin{lstlisting}
from firedrake import *
msh = UnitSquareMesh(16, 16)
x, y = SpatialCoordinate(msh)
f = cos(pi * x) * sin(3 * pi * y)
V = FunctionSpace(msh, 'CG', 2)
u = Function(V)
v = TestFunction(V)
F = inner(grad(u), grad(v)) * dx - inner(f, v) * dx
bcs = Dirichlet(V, 0, 'on_boundary')
solve(F==0, u, bcs=bcs)
  \end{lstlisting}
  \caption{Firedrake listing for the Poisson equation.}
  \label{fig:fdpoisson}
  \Description{Code listing for the Poisson equation in Firedrake.}
\end{figure}

UFL, however, has historically lacked an abstraction to directly describe time-dependent problems and their evolution.
Hence, users have typically hand-written their own time stepping loops using simple methods like backward Euler or Crank-Nicolson, although Dune has made some progress at enabling Galerkin-in-time methods in~\cite{versbach2023theoretical}.
Changing the time stepping scheme then becomes invasive and requires rewriting a great deal of the application code.
We contrast this with the ease with which one may change the order of spatial approximation by varying an argument to the \lstinline{FunctionSpace} constructor in Figure~\ref{fig:fdpoisson}.

There has, moreover, been a renewed interest in higher-order, fully implicit \RK{} methods for integrating spatial finite element discretizations in time.
These methods were long dismissed as impractical, as they require a rather complex implementation path and lead to large, coupled algebraic systems for all stages within a single time step.
On the other hand, these methods lack the order barriers present in explicit and diagonally implicit methods.
Collocation-type methods can provide superconvergent algebraically stable methods of any order.  They also have high  \emph{stage order}, so they do not suffer from the same level of order reduction as diagonally implicit methods~\cite{frank1985order,wanner1996solving}.
Recent work in~\cite{masud2021new,pazner2017stage,southworth2022fast1,southworth2022fast2,vanlent2005,abu2022monolithic} has shown that effective algebraic solvers can make fully implicit methods competitive or even superior to other approaches in practice.

\Irksome~\cite{farrell2021irksome,kirby2024extending} addresses this range of concerns within the Firedrake project.
First, it provides an abstraction for time-dependent problems, enabling users to write a semidiscrete variational form rather than explicitly defining a time stepping scheme, as shown in Figure~\ref{fig:fdirkheat} for the heat equation.
Here, the UFL language has been extended with the \lstinline{Dt} operator representing time differentiation, and we have new classes for time steppers and Butcher tableaux.
Second, although implicit \RK{} methods require the solution of a complex, stage-coupled variational problem, \Irksome{} automatically generates this from the semidiscrete problem description.
Finally, Firedrake possesses a rich interface to PETSc~\cite{kirby2018solver}.
It is possible to extend the large suite of solvers and preconditioners PETSc provides with new operations defined with Firedrake.
For example, patch-based additive Schwarz methods are included in~\cite{farrell2021pcpatch}, and 
we apply a facility for preconditioning with auxiliary operators to generate the novel \RK{} preconditioners of~\citet{masud2021new} in our recent work~\cite{kirby2024extending}.

\begin{figure}
  \begin{lstlisting}
from firedrake import *
from irksome import Dt, TimeStepper, RadauIIA
msh = UnitSquareMesh(16, 16)
x, y = SpatialCoordinate(msh)
t = Constant(0)
f = cos(pi * x) * sin(3 * pi * y) * cos(t)
V = FunctionSpace(msh, 'CG', 2)
u = Function(V)
v = TestFunction(V)
F = (inner(Dt(u), v) * dx + inner(grad(u), grad(v)) * dx
     - inner(f, v) * dx)
bcs = Dirichlet(V, 0, 'on_boundary')
butcher_tableau = RadauIIA(2)
dt = Constant(1/16)
stepper = TimeStepper(F, butcher_tableau, u, t, dt, bcs=bcs)
while float(t) < 1.0:
    stepper.advance()
    t.assign(float(t) + float(dt))
  \end{lstlisting}
  \caption{Firedrake/\Irksome~listing for the heat equation}
  \label{fig:fdirkheat}
  \Description{Code listing for the heat equation in Irksome.}
\end{figure}

In this paper, we build on the recent successes of \Irksome{} and introduce support for \RKN{} methods~\cite{nystrom1925, wanner1993solving} for problems with second-order time derivatives.
We have now enabled the user code in Figure~\ref{fig:fdnrkwave} for the wave equation, extending \lstinline{Dt} to support higher-order derivatives and adding Nystr\"om time steppers.
While such problems can be programmatically reduced to first-order systems and integrated with a standard \RK{} method, directly tackling the second-order system allows a UFL problem description more closely aligned with the mathematics.
Avoiding the auxiliary variables saves memory and leads to simpler algebraic systems.  As we will show, this leads to performance improvements when solving in second-order form.

\begin{figure}
  \begin{lstlisting}
from firedrake import *
from irksome import Dt, NystromTimeStepper, GaussLegendre
msh = UnitSquareMesh(16, 16)
x, y = SpatialCoordinate(msh)
t = Constant(0)
f = cos(pi * x) * sin(3 * pi * y) * cos(t)
V = FunctionSpace(msh, 'CG', 2)
u = Function(V)
v = TestFunction(V)
F = (inner(Dt(u, 2), v) * dx + inner(grad(u), grad(v)) * dx
     - inner(f, v) * dx)
bcs = Dirichlet(V, 0, 'on_boundary')
butcher_tableau = GaussLegendre(2)
dt = Constant(1/16)
stepper = NystromTimeStepper(F, butcher_tableau, u, t, dt, bcs=bcs)
while float(t) < 1.0:
    stepper.advance()
    t.assign(float(t) + float(dt))
  \end{lstlisting}
  \caption{Firedrake/\Irksome~listing for the wave equation}
  \label{fig:fdnrkwave}
  \Description{Code listing for the wave equation in Irksome.}  
\end{figure}

In the rest of the paper, we introduce \RK{} and \RKN{} methods in Section~\ref{sec:rk}.
In addition, we describe the treatment of essential boundary conditions when these methods are applied to finite element discretizations of PDEs.
Then, in Section~\ref{sec:models}, we describe a handful of particular models.
These include the standard wave equation, a biharmonic wave equation, and a fully dynamic poroelastic model.
As with \RK{} methods, a common challenge for implicit \RKN{} methods is the solution of a large, fully-coupled algebraic system for the stages, and we describe two distinct approaches in Section~\ref{sec:alg}.
A description of how \RKN{} methods and solvers are implemented in Irksome is given in Section~\ref{sec:imp}, and
we evaluate our implementation on the given model problems in Section~\ref{sec:numres} before concluding in Section~\ref{sec:conc}.

\section{\RK{} and \RKN{} methods}
\label{sec:rk}
Given domain $\Omega$, we are interested in the solution of (systems of) time-dependent PDEs of the form
\begin{equation}\label{eq:generic_PDE}
  u_{tt} = F(t, x, u, u_t) \text{ in }\Omega,
\end{equation}
where $F$ can include spatial derivatives of $u$ and $u_t$, but not further temporal derivatives,
combined with appropriate initial conditions, $u(0,x) = g(x)$ and $u_t(0,x) = \gamma(x)$ for all $x\in\Omega$.  Depending on the form of $F$ (in particular its dependence on spatial derivatives of $u$), this form will be accompanied by one or more boundary conditions, specifying values of $u$ and/or its derivatives along segments $\Gamma \subset \partial\Omega$, such as $u(t,x) = h(t,x)$ for $x\in\Gamma$.  We note that our methods are able to solve much more complex systems of PDEs than those in~\eqref{eq:generic_PDE}, including differential-algebraic equations (DAEs) and coupled systems mixing first- and second-order time derivatives, but we focus on the simplest case for exposition in this section.

We consider the spatial discretization of~\eqref{eq:generic_PDE} using finite element methods, with specific examples given below in~\cref{sec:models}.  Given some finite element space, $V_h$, with basis $\{\psi_\ell\}_{\ell=1}^{\text{dim} V_h}$, we write $\displaystyle u(t,x) = \sum_{\ell=1}^{\text{dim} V_h} y_\ell(t)\psi_\ell(x)$.  We then rewrite~\eqref{eq:generic_PDE} in discretized variational form, as
\begin{equation}\label{eq:generic_variational}
  (u_{tt},v) = F(t, x, u, u_t; v) \text{ for all }v\in V_h,
\end{equation}
where $(u,v) = \int_\Omega uv$ is the standard $L^2(\Omega)$ inner product.  As is usual in finite element methods, we then generally integrate-by-parts on the right-hand side (accounting for boundary conditions) and substitute in the basis functions to get a set of ordinary differential equations for the coefficients, $\{y_\ell(t)\}_{\ell=1}^{\text{dim} V_h}$.  While the mass matrix that arises in the left-hand side is important in the details of the approaches that follow, we can always invert it to rewrite the resulting semi-discrete system of equations as $\bfy_{tt}(t) = \bff(t, \bfy, \bfy_t)$ with initial conditions given by interpolating (or projecting) those on $u$ to appropriate initial conditions on $\bfy$.  We consider two approaches for integrating the resulting system of ODEs, either rewriting it as a first-order system of equations and integrating using standard \RK{} methods or directly integrating the second-order form using \RKN{} methods.

\subsection{Methods for first-order systems}
Thus, we begin by considering systems of first-order ordinary differential equations of the form
\begin{equation}
  \label{eq:firstorderode}
  \begin{split}
    \bfy_t(t) & = \bff(t, \bfy), \\
    \bfy(0) & = \bfy^0.
  \end{split}
\end{equation}
where $\bff: \left(0, T\right] \times \mathbb{R}^m \rightarrow \mathbb{R}^m$ and the solution $\bfy$ is a function mapping $(0, T] \rightarrow \mathbb{R}^m$.

    \RK{} methods comprise a vast family of methods for approximating solutions to~\eqref{eq:firstorderode}.
    We partition $[0, T]$ by $0 = t^0 < t^1 < \dots < t^N \equiv T$ and put
    $\Delta t^n = t^n - t^{n-1}$.  In the sequel, we shall assume uniform time steps for simplicity of notation, although (unlike multi-step methods) there is no difficulty in formulating \RK{} methods with variable steps. 
    Whatever the exact properties of a particular method are, these methods are uniformly expressed via a Butcher tableau
\begin{equation}
  \label{eq:butcher_tableau}
  \begin{array}{c|c}
    \bfc & A \\ \hline
    & \bfb
  \end{array},
\end{equation}
where the vectors $\bfb,\bfc \in \mathbb{R}^s$ contain the entries $b_i$ and
$c_i$, respectively, and $A \in \mathbb{R}^{s \times s}$ contains the
entries $a_{ij}$.  The natural number $s$ refers to the number of \emph{stages} in the RK method.

Given some approximation $\bfy^n \approx \bfy(t^{n})$, a \RK{} method seeks
$\bfy^{n+1} \approx \bfy(t^{n+1})$ given by
\begin{equation}
  \bfy^{n+1} = \bfy^n + \Delta t \sum_{i=1}^s b_i \bfk^{(i)},
\end{equation}
where the $\bfk^{(i)} \in \mathbb{R}^m$ for $1 \leq i \leq s$ are defined by the algebraic system
\begin{equation}
  \label{eq:stagesystem}
  \bfk^{(i)} = \bff\left(t^n + c_i \Delta t, \bfy^n + \Delta t \sum_{j=1}^s a_{ij} \bfk^{(j)} \right).
\end{equation}

\emph{Explicit} methods correspond to the case where $A$ is strictly lower triangular.  That is, $a_{ij} = 0$ for $1 \leq i \leq j \leq s$.  Then,~\eqref{eq:stagesystem} becomes
\begin{equation}
  \label{eq:expstagesystem}
  \bfk^{(i)} = \bff\left(t^n + c_i \Delta t, \bfy^n + \Delta t \sum_{j=1}^{i-1} a_{ij} \bfk^{(j)} \right),
\end{equation}
and each $\bfk^{(i)}$ can be computed by evaluating $\bff$ in terms of previously computed stages.

When $A$ is lower triangular, the method is called \emph{diagonally implicit}, and
\begin{equation}
  \label{eq:dirkstagesystem}
  \bfk^{(i)} = \bff\left(t^n + c_i \Delta t, \bfy^n + \Delta t \sum_{j=1}^{i} a_{ij} \bfk^{(j)} \right).
\end{equation}
That is, each $\bfk^{(i)}$ is computed by solving an algebraic system, with the solution of that system depending on previously computed $\bfk^{(j)}$ for $1 \leq j < i$.
Moreover, this algebraic system is quite similar to that required by the simpler backward Euler method.  Consequently, an existing implementation
of backward Euler can be readily repurposed to compute diagonally implicit methods.  This feature, shared by implicit multistep methods, may explain the common
application of such methods when higher temporal accuracy is desired.
However, implicit multistep methods suffer from the second Dahlquist barrier~\cite{dahlquist1963special} -- A-stable methods cannot exceed second order.
More subtly, diagonally implicit methods suffer from \emph{order reduction} in the presence of stiffness.
Here, the accuracy of the methods reduces to the stage order (the truncation error of the individual stages), and this cannot exceed two for
irreducible DIRK methods~\cite{wanner1996solving}. 

There exist \emph{fully implicit} \RK{} methods that do not suffer from these order barriers. %
Excellent theoretical properties of fully implicit methods, especially those arising from polynomial collocation, have been well known for some time~\cite{wanner1996solving}.
For example, methods based on Gauss-Legendre collocation are known to be algebraically stable, with an $s$-stage method offering formal accuracy of order $2s$.
Moreover, these methods are symplectic for all $s$.
The RadauIIA methods are algebraically stable with the $s$-stage method offering formal accuracy of order $2s-1$.
Their $L$-stability makes them well-suited for strongly stiff problems.
The $s$-stage instances of each family have stage order $s$, so the effects of order reduction in the presence of stiffness is far less pronounced than
for diagonally implicit methods, and they can give superior accuracy to DIRKs with high weak stage order~\cite{biswas2023design} for some challenging benchmark problems~\cite{kirby2024extending}.

However, for fully implicit methods,~\eqref{eq:stagesystem} requires solving a single algebraic system coupling all the stages simultaneously.
This system is far more complicated than the per-stage systems for diagonally implicit methods.
Especially in the context of discretized PDEs, there has been much recent progress in solving that stage-coupled system.
One class of methods proposes preconditioners that are, in effect, diagonally implicit and so allow reuse of effective backward Euler solvers~\cite{masud2021new, nilssen2011order}.
Alternatively, monolithic multigrid methods, with relaxation schemes that couple the stages together, have also proven very powerful~\cite{vanlent2005, farrell2021irksome, abu2022monolithic}.
In fact, the numerical experiments in papers such as~\cite{kirby2024extending} show that these modern techniques can make fully implicit methods competitive with the simpler alternatives.  Indeed, fully implicit methods can even give greater accuracy for a given run-time.

\subsection{\RKN{} methods}
This paper is primarily focused on second-order differential equations of the form
\begin{equation}
  \label{eq:ode}
  \begin{split}
    \bfy_{tt}(t) & = \bff(t, \bfy, \bfy_t), \\
    \bfy(0) & = \bfy^0, \\
    \bfy_t(0) & = \bfy_t^0,
    \end{split}
\end{equation}
where $\bff:\left(0, T\right] \times \mathbb{R}^m \times \mathbb{R}^m \rightarrow \mathbb{R}^m$ and the solution $\bfy : \left(0, T\right] \rightarrow \mathbb{R}^m$.

Introducing new variables $\bfy^{(1)} = \bfy$ and $\bfy^{(2)} = \bfy_t$ enables us to rewrite this equation as a first order system in the standard way:
\begin{equation}
  \begin{split}
    \bfy^{(1)}_t & = \bfy^{(2)}, \\
    \bfy^{(2)}_t & = \bff(t, \bfy^{(1)}, \bfy^{(2)}).
  \end{split}
  \label{eq:odesystem}
\end{equation}

Applying a generic \RK{} method with tableau~\eqref{eq:butcher_tableau} to this rewritten system requires two sets of stages that we label $\{ \bfk^{(1, i)} \}_{i=1}^s$ and $\{ \bfk^{(2, i)} \}_{i=1}^s$, defined by the algebraic system
\begin{equation}
  \begin{split}
    \bfk^{(1, i)} & = \bfy^{(2,n)} + \Delta t \sum_{j=1}^s a_{ij} \bfk^{(2, j)}, \\
    \bfk^{(2, i)} & = \bff\left(t^n + c_i \Delta t,
    \bfy^{(1,n)} + \Delta t \sum_{j=1}^s a_{ij} \bfk^{(1, j)},
    \bfy^{(2,n)} + \Delta t \sum_{j=1}^s a_{ij} \bfk^{(2, j)} \right),
  \end{split}
\label{eq:fosysstage}
\end{equation}
and we update the solution by
\begin{equation}
  \begin{split}
    \bfy^{(1,n+1)} & = \bfy^{(1,n)} + \Delta t \sum_{i=1}^s b_i \bfk^{(1, i)}, \\
    \bfy^{(2,n+1)} & = \bfy^{(2,n)} + \Delta t \sum_{i=1}^s b_i \bfk^{(2, i)}.
  \end{split}
  \label{eq:fosysupdate}
\end{equation}

We can simplify this process by substituting the first equation in~\eqref{eq:fosysstage} into the remaining equations of~\eqref{eq:fosysstage} and~\eqref{eq:fosysupdate}.  Defining $\bfkappa^{(i)} \equiv \bfk^{(2, i)}$, we arrive at the much smaller algebraic system
\begin{equation}
  \label{eq:rknstageeq}
    \bfkappa^{(i)} = \bff\left(t^n + c_i \Delta t,
    \bfy^{(1,n)} + c_i \Delta t \bfy^{(2,n)} + \left( \Delta t \right)^2 \sum_{j=1}^s \overline{a}_{ij} \bfkappa^{(j)},
    \bfy^{(2,n)} + \Delta t \sum_{j=1}^s a_{ij} \bfkappa^{(j)} \right),
\end{equation}
where
\begin{equation}
  \label{eq:abar}
  \overline{a}_{ij} = \sum_{\ell=1}^s a_{i\ell} a_{\ell j}.
\end{equation}
The update formulae simplify to
\begin{equation}
  \label{eq:rknupdate}
  \begin{split}
    \bfy^{(1,n+1)} & = \bfy^{(1,n)} + \Delta t \bfy^{(2,n)} + \left( \Delta t \right)^2 \sum_{i=1}^s \overline{b}_i \bfkappa^{(i)}, \\
    \bfy^{(2,n+1)} & = \bfy^{(2,n)} + \Delta t \sum_{i=1}^s b_i \bfkappa^{(i)},
  \end{split}
\end{equation}
where
\begin{equation}
  \label{eq:bbar}
  \overline{b}_i = \sum_{j=1}^s a_{ji} b_j.
\end{equation}
This formulation has certain advantages over applying a \RK{} method to the first order system.  Most notably, it only requires one set of stage variables rather than two.  This reduces the memory requirement and, for implicit methods, greatly simplifies the relevant algebraic system to solve.  We note that \RK{} methods can alternately be formulated in terms of approximations to the solution at the stage times.  However, this "stage-value" formulation is less amenable to handling general second-order equations, as we sketch in the appendix, so we do not pursue this here.

With the definitions~\eqref{eq:abar},~\eqref{eq:bbar}, we can introduce notation for an extended tableau
\begin{equation}
  \label{eq:rkn_butcher_tableau}
  \begin{array}{c|c|c}
    \bfc & \overline{A} & A \\ \hline
    & \overline{\bfb} & \bfb
  \end{array}.
\end{equation}
Equations~\eqref{eq:abar} and~\eqref{eq:bbar} tell us that when our method is derived from applying an existing RK method to the first-order system,
\begin{equation}\label{eq:inherited_nystrom}
  \begin{split}
    \overline{A} & = A^2, \\
    \overline{\bfb} & = A^T \bfb.
    \end{split}
\end{equation}
For example, the extended tableau for the two-stage Gauss--Legendre method is
\begin{equation}
  \label{eq:gl2ext}
  \begin{array}{c|cc|cc}
    \tfrac{1}{2} - \tfrac{\sqrt{3}}{6} & \tfrac{1}{24} & \tfrac{1}{8} - \tfrac{\sqrt{3}}{12} & \tfrac{1}{4} & \tfrac{1}{4} - \tfrac{\sqrt{3}}{6} \\
    \tfrac{1}{2} + \tfrac{\sqrt{3}}{6} & \tfrac{1}{8} + \tfrac{\sqrt{3}}{12} & \tfrac{1}{24} & \tfrac{1}{4} + \tfrac{\sqrt{3}}{6} & \tfrac{1}{4} \\ \hline
    & \tfrac{1}{4} + \tfrac{\sqrt{3}}{12} & \tfrac{1}{4} - \tfrac{\sqrt{3}}{12}
    & \tfrac{1}{2} & \tfrac{1}{2}
  \end{array}.
\end{equation}

However, methods of the form~\eqref{eq:rknstageeq},~\eqref{eq:rknupdate} exist that are not derivable from RK methods for a first-order system, so the formalism of an extended tableau~\eqref{eq:rkn_butcher_tableau} is appropriate.  Nystr\"om proposed such a method in~\cite{nystrom1925}:

\begin{equation}\label{eq:nystrom_scheme}
    \begin{array}{c|cccc|cccc}
      0 & & & & & & & & \\
      \tfrac{1}{2} & \tfrac{1}{8} & & & & \tfrac{1}{2} & & & \\
      \tfrac{1}{2} & \tfrac{1}{8} & 0 & & & 0 & \tfrac{1}{2} &  & \\
      1 & 0 & 0 & \tfrac{1}{2} & & 0 & 0 & 1 & \\ \hline
      & \tfrac{1}{6} & \tfrac{1}{6} &  \tfrac{1}{6} & 0 &  \tfrac{1}{6} &  \tfrac{1}{3} & \tfrac{1}{3} &  \tfrac{1}{6} \\
  \end{array}.
\end{equation}

\subsection{Boundary conditions}

While, for the most part, we can directly apply \RK{} and \RKN{} methods to the semi-discretized form in~\eqref{eq:generic_variational} (or that obtained after integrating by parts), the strong imposition of time-dependent boundary conditions requires care.
Accurately imposing Dirichlet boundary conditions on \RK{} discretizations of time-dependent PDEs is recognized to be a difficult task~\cite{carpenter1995theoretical, pathria1997correct}, which is not made easier by the additional derivatives appearing in second-order equations.  Here, we outline three consistent approaches to strongly impose boundary conditions, given time-dependent Dirichlet boundary conditions prescribing $u(t,x) = h(t,x)$ for $x\in \Gamma \subset \partial\Omega$.  In all that follows, we assume that the given boundary condition data, $h(t,x)$, has been interpolated (or projected) into the finite element space, so that we can interchangeably think about $h(t,x)$ itself and its coefficients for the subset of the finite element basis functions, $\psi_\ell(x)$, that are nonzero on $\Gamma \subset \partial\Omega$.

Our first approach, denoted as \lstinline{ODE} boundary conditions, makes use of the direct analogy between the \RKN{} stages, $\bfkappa^{(i)}$, and the second derivative $u_{tt}(t,\cdot)$  In particular, for coefficients of finite element basis functions that are nonzero on $\Gamma$, the prescription of $u(t,x) = h(t,x)$ naturally also prescribes $u_{tt}(t,x) = h_{tt}(t,x)$, so we can strongly enforce a boundary condition on stage $i$ by requiring that the associated basis coefficient in $\bfkappa^{(i)}$ match that for $h_{tt}(t^n+c_i\Delta t,x)$.  We note two shortcomings of this approach.  First, it offers no guarantee that $u(t^{n+1},x)$ will interpolate $h(t^{n+1},x)$ on $\Gamma$, since these boundary conditions are imposed without any accommodation for possible discontinuities (in time) of $h(t,x)$~\cite{kirby2024extending}.  Indeed, it offers no guarantee that the boundary conditions on $\bfkappa^{(i)}$ are even well-defined, unless $h_{tt}(t,x)$ is assumed to be continuous in time.  Secondly, while the ODE enforced on the boundary values is consistent with the Dirichlet condition, it is inconsistent with the true nature of Dirichlet boundary conditions, which impose a differential-algebraic equation (DAE) structure on the spatially semi-discretized form.

Our second approach, denoted as \lstinline{DAE} boundary conditions, arises from directly interpreting the ODE plus boundary conditions as a DAE.  In this setting, we impose the values of the stage approximation,
\begin{equation}
  \left(\bfy^{(1,n)} + c_i\Delta t \bfy^{(2,n)} + \left(\Delta t\right)^2\sum_{j=1}^s \bar{a}_{ij}\bfkappa^{(j)}\right)_\ell,
\end{equation}
for each basis function, $\psi_\ell(x)$, that is nonzero on the boundary, prescribing that it match the corresponding coefficient in the basis expansion of $h(t^n + c_i\Delta t,x)$ for $1\leq i \leq s$.
This gives a coupled system of $s$ equations for the values $\{\kappa^{(j)}_\ell\}_{j=1}^s$ for each such basis function, $\psi_\ell$, that we aim to solve to properly impose the boundary conditions in DAE form.  Clearly our ability to do this relies on the invertibility of $\bar{A}$.  When the \RKN{} tableau is inherited from an implicit \RK{} method, this is generally the case; however, when the method is explicit, we know that $\bar{A}$ will not be invertible.  For the explicit \RK{} case, we can generally still apply boundary conditions using the DAE viewpoint, by ``shifting'' the set of stage values at which we enforce the boundary condition, since there is no mechanism for enforcing a boundary condition on the explicit first stage.  Instead, we can impose the boundary condition as a DAE by requiring it to hold for stages $2$ through $s$ as well as for the approximation at time $t^{n+1}$.  When the \RK{} tableau has rank $s-1$ and the stage reconstruction is linearly independent, this is well-posed (and, indeed, this is what is implemented for DAE boundary conditions on explicit \RK{} schemes in \Irksome.

For explicit \RKN{} methods, however, it is possible that $\bar{A}$ will have rank $s-2$ and, consequently, the approach from the explicit \RK{} case cannot be applied.  This is the case, for example, in Nystr\"om's scheme from~\eqref{eq:nystrom_scheme}, where $\bar{A}$ clearly has rank 2 for a 4-stage scheme.  Similarly, if we derive the \RKN{} scheme from an explicit \RK{} scheme, as in~\eqref{eq:inherited_nystrom}, we expect $A$ to be of rank $s-1$, and $\bar{A} = A^2$ should have rank $s-2$, due to the strictly lower-triangular structure of $A$.  To overcome this, we first revisit the equivalent first-order system in~\eqref{eq:odesystem}, and note that the \RKN{} stages $\bfkappa^{(i)}$ correspond to the stage-derivative approximations for $\bfy_2 = \bfy_1'$.  Thus, when we seek to impose boundary conditions on $\bfkappa^{(i)}$ consistent with the DAE framework, it is perhaps more natural to impose them on $\bfy_2$ in the equivalent first-order system.  This means that we apply boundary conditions on $u_t(t,x) = h_t(t,x)$ for $x\in\Gamma\subset\partial\Omega$.  For the implicit case, we then use the stage approximations and set
\[
\left(\bfy^{(2,n)} + \Delta t\sum_{j=1}^s {a}_{ij}\bfkappa^{(j)}\right)_\ell
\]
to match the basis expansion of $h_t(t^n + c_i\Delta t, x)$ for each basis function, $\ell$, that is nonzero on the boundary, for $1 \leq i \leq s$.  As this again relies on the invertibility of $A$ to solve for $\{\kappa^{(j)}_\ell\}_{j=1}^s$, for the explicit case, we use the same approach as described above for explicit \RK{} schemes, imposing values on the stage approximations for $2\leq i \leq s$ and a similar requirement for the approximation at time $t^{n+1}$.  We denote this type of boundary condition as \lstinline{dDAE} boundary conditions (for differentiated DAE boundary conditions).

\section{Some model partial differential equations}\label{sec:models}
Here, we present a suite of semidiscrete finite element methods for some model partial differential equations with second-order temporal derivatives.
A prototypical problem is the wave equation, posed on some domain $\Omega \subset \mathbb{R}^d$ with $1 \leq d \leq 3$ and time interval $[0, T]$.  To fix ideas, we take homogeneous Dirichlet boundary conditions.
After decomposing $\Omega$ into some triangulation and defining a finite-dimensional space $V_h \subset H_0^1(\Omega)$ comprising continuous piecewise polynomials of some degree $k \geq 1$ that also vanish on the boundary,
we have a time-dependent variational problem
\begin{equation}
  \label{eq:semidwave}
  \left( u_{tt}, v \right) + \left( \nabla u , \nabla v \right) = 0,
\end{equation}
where $\left( \cdot, \cdot \right)$ denotes the standard $L^2$ inner product.
The system is completed with initial conditions on $u$ and its first temporal derivative:
\begin{equation}
  \begin{split}
    u(0, \cdot) & = g(\cdot), \\
    u_t(0, \cdot) & = \gamma(\cdot).
  \end{split}
\end{equation}

By expanding the solution
\begin{equation}
  u(t, \cdot) = \sum_{\ell=1}^{\dim V_h} y_\ell(t) \psi_\ell
\end{equation}
in some basis $\{ \psi_\ell \}_{\ell=1}^{\dim V_h}$ for $V_h$, we arrive at a system of ODEs
\begin{equation}
  \label{eq:waveODE}
  M \bfy_{tt} = - K \bfy,
\end{equation}
where $M_{k\ell} = \left( \psi_\ell, \psi_k \right)$ and
$K_{k\ell} = \left( \nabla \psi_\ell , \nabla \psi_k \right)$ are the standard mass and stiffness matrices.  One also requires approximations to the initial data.

Absent forcing, choosing $v = u_t$ in~\eqref{eq:semidwave} leads to the standard energy conservation result
\begin{equation}
  \label{eq:waveconserv}
  \frac{d}{dt} \left( \frac{1}{2} \left\| u_t \right\|^2 + \frac{1}{2} \left\| \nabla u \right\|^2 \right) = 0,
\end{equation}
and this can be written equivalently as
\begin{equation}
  \frac{d}{dt} \left( \frac{1}{2} \bfy_t^\top M \bfy_t + \frac{1}{2} \bfy^\top K \bfy \right) = 0.
\end{equation}

The Gauss-Legendre schemes exactly maintain this energy conservation.  The lowest-order case is equivalent to the implicit midpoint rule, but higher-order methods require the solution of a stage-coupled algebraic system at each time step.
To illustrate this, we suppose that we have some $u^n, u_t^n \in V_h$ approximating the solution and its temporal derivative at time level $n$.
Applying a two-stage \RKN{} method such as that of~\eqref{eq:gl2ext} to the wave equation gives the variational problem of seeking $\kappa_1, \kappa_2 \in V_h$ such that 
\begin{equation}
  \label{eq:2stagerknwave}
  \begin{split}
    \left( \kappa_1, v_1 \right)
    + \left( \nabla
    \left( u^n + c_1 \left( \Delta t \right) u_t^n
    + \left( \Delta t \right)^2 \left(
    \overline{a}_{11} \kappa_1
    +\overline{a}_{12} \kappa_2  \right)
    \right)
    , \nabla v_1 \right) & = 0, \\
    \left( \kappa_2, v_2 \right)
    + \left( \nabla
    \left( u^n + c_2 \left( \Delta t \right) u_t^n
    + \left( \Delta t \right)^2 \left(
    \overline{a}_{21} \kappa_1
    +\overline{a}_{22} \kappa_2  \right)
    \right)
    , \nabla v_2 \right) & = 0,
  \end{split}
\end{equation}
for all $v_1, v_2 \in V_h$.   After solving this system we compute
$u^{n+1}$ and $u_t^{n+1}$ by~\eqref{eq:rknupdate}.

By expanding the stages $\kappa_i = \sum_{\ell=1}^{\dim V_h} k^i_\ell \psi_\ell$ in the basis for $V_h$, we arrive at a block-structured linear system
\begin{equation}
  \begin{bmatrix} M + \overline{a}_{11} \Delta t^2 K & \overline{a}_{12} \Delta t^2 K \\
    \overline{a}_{21} \Delta t^2 K & M + \overline{a}_{22}\Delta t^2 K
  \end{bmatrix}
  \begin{bmatrix} \bfk^{(1)} \\ \bfk^{(2)}
  \end{bmatrix}
  = \begin{bmatrix} \bfF^{(1)} \\ \bfF^{(2)}
  \end{bmatrix},
  \label{eq:blocksys}
\end{equation}
where $\bfk^{(1)}$ and $\bfk^{(2)}$ are the vectors of coefficients in the basis expansion, and $\bfF^{(1)}$ and $\bfF^{(2)}$ gather the known terms from~\eqref{eq:2stagerknwave}, expressed in the same spatial basis.
We note that the variational problem and associated matrix equation depends strongly on the \emph{number} of stages in the \RKN{} method, although the structure for a given number of stages varies very little with the particular values in the tableau (unless of course, one obtains an explicit or diagonally implicit method).
We also note that the system has only $2 \times \dim V_h$ unknowns, rather than the $4 \times \dim V_h$ unknowns required if we were to rewrite the system in first-order form.

Explicit dependence of the PDE on the time derivative $u_t$ leads to a more complicated variational problem for the stages.
Although our goal is to apply \RKN{} methods to more complex problems like fully dynamic poroelasticity, it helps to see how first derivatives enter into a much simpler equation.
So, consider next the telegraph equation, written in variational form as
\begin{equation}
  \left( u_{tt}, v \right) + \left( u_t , v \right) + \left( \nabla u , \nabla v \right) = 0.
\end{equation}
The energy functional in~\eqref{eq:waveconserv}, conserved for the wave equation, can be shown to decay exponentially over time for the telegraph equation, although one does not gain smoothness like with the heat equation.
Energy decay makes B-stability an attractive feature for implicit methods.
Hence, Gauss-Legendre methods would still be suitable, although any advantage afforded symplecticity would be unclear here.

Applying a two-stage \RKN{} method to the telegraph equation gives a coupled variational problem for $\kappa_1, \kappa_2 \in V_h$ such that 
\begin{equation}
  \label{eq:2stagerkntele}
  \begin{split}
    \left( \kappa_1, v_1 \right)
    + \left( u_1^n + \Delta t \left( a_{11} \kappa_1 + a_{12} \kappa_2 \right), v_1 \right)
    + \left( \nabla
    \left( u_1^n + c_1 \left( \Delta t \right) u_2^n
    + \left( \Delta t \right)^2 \left(
    \overline{a}_{11} \kappa_1
    +\overline{a}_{12} \kappa_2  \right)
    \right)
    , \nabla v_1 \right) & = 0, \\
    \left( \kappa_2, v_2 \right)
    + \left( u_1^n + \Delta t \left( a_{21} \kappa_1 + a_{22} \kappa_2 \right), v_2 \right)
    + \left( \nabla
    \left( u_1^n + c_2 \left( \Delta t \right) u_2^n
    + \left( \Delta t \right)^2 \left(
    \overline{a}_{21} \kappa_1
    +\overline{a}_{22} \kappa_2  \right)
    \right)
    , \nabla v_2 \right) & = 0,
  \end{split}
\end{equation}
and this gives rise to a block-structured system
\begin{equation}
  \begin{bmatrix} \left( 1 + a_{11} \Delta t\right) M + \overline{a}_{11} \Delta t^2 K & a_{12} \Delta t M + \overline{a}_{12} \Delta t^2 K \\
    a_{21} \Delta t M + \overline{a}_{21} \Delta t^2 K & \left( 1 + a_{22} \Delta t \right) M + \overline{a}_{22}\Delta t^2 K
  \end{bmatrix}
  \begin{bmatrix} \bfk^{(1)} \\ \bfk^{(2)}
  \end{bmatrix}
  = \begin{bmatrix} \bfF^{(1)} \\ \bfF^{(2)}
  \end{bmatrix},
\end{equation}
again with $\bfF^{(1)}$ and $\bfF^{(2)}$ gathering known terms expressed in the finite element basis.


We also consider a second-order equation with fourth-order spatial derivatives -- a dynamic plate problem
\begin{equation}
  \label{eq:dynplate}
  u_{tt} + \Delta^2 u = 0
\end{equation}
on some domain $\Omega$, subject to clamped boundary conditions
\begin{equation}
  u|_\Omega = \frac{\partial u}{\partial n}|_{\Omega} = 0.
\end{equation}
and initial conditions on the solution $u$ and its time derivative $u_t$ on the entire domain.
This equation models the dynamic evolution of a thin plate.  Standard arguments show the conservation property
\begin{equation}
  \label{eq:dynplateenergy}
  \frac{d}{dt} \left[ \frac{1}{2} \| u_t \|^2 + \frac{1}{2} \| \Delta u \|^2 \right] = 0.
\end{equation}

Here, the natural space for a variational formulation is
$H^2_0(\Omega)$, the space of $H^2$ functions with vanishing trace.
We take some finite-dimensional $V_h \subset H^2_0(\Omega)$ and seek a discrete solution $u: [0, T] \rightarrow V_h$ such that
\begin{equation}
  \left( u_{tt}, v \right) + \left( \Delta u , \Delta v \right) = 0
\end{equation}
for all $v \in V_h$.
We obtain $V_h$ by the Hsieh-Clough-Tocher macroelement~\cite{clough1965finite,brubeck2025fiat} on a triangulation of $\Omega$.  

Although we have wave-like behavior and a conserved quantity, explicit methods will be far less competitive for this equation than the standard wave equation.
The higher-order spatial derivatives lead to a time step restriction for stability of $\Delta t = \mathcal{O}(h^2)$, rather than $\mathcal{O}(h)$ for the wave equation, greatly increasing the cost of explicit time stepping to any fixed $\mathcal{O}(1)$ final time.
In contrast, Gauss-Legendre methods provide an unconditionally stable method that maintains energy conservation in the discrete analog of~\eqref{eq:dynplateenergy} without this severe time step restriction.

As a third example, we consider a fully dynamic Biot model of poroelasticity~\cite{10.1121/1.1908239, BOTH2022114183, kraus2024analysis}.
This is a system of PDEs modeling the interaction between fluid flow and elastic deformation of a porous medium.
It is of mixed hyperbolic/parabolic type, with hyperbolic elastodynamic deformation coupled to parabolic Darcy flow and mass transport.
We solve the system
\begin{equation}\label{eq:biot_system}
  \begin{split}
    \overline{\rho} \boldsymbol{u}_{tt}
    - 2 \mu \nabla \cdot \epsilon \left( \boldsymbol{u} \right)
    - \lambda \nabla \left( \nabla \cdot \boldsymbol{u} \right)
    + \alpha \nabla p
    + \rho_f \boldsymbol{w}_t & = \boldsymbol{f}, \\
    \rho_f \boldsymbol{u}_{tt}
    + \rho_w \boldsymbol{w}_t
    + K^{-1} \boldsymbol{w} + \nabla p & = \boldsymbol{g}, \\
    s_0 p_t
    + \alpha \nabla \cdot \boldsymbol{u}_t
    + \nabla \cdot \boldsymbol{w} & = 0,
  \end{split}
\end{equation}
for the deformation $\boldsymbol{u}$, fluid velocity  $\boldsymbol{w}$, and fluid pressure $p$.  The physical parameters of the problem depend on the medium's porosity, $\phi_0$, which may vary spatially.
The total density $\overline{\rho} = (1-\phi_0)\rho_s + \phi_0 \rho_f$, where $\rho_s$ and $\rho_f$ are the density of the solid and fluid, respectively, while $\rho_w \geq \phi_0^{-1}\rho_f$ is the effective fluid density.
The Biot-Willis parameter $\alpha$ lies between $\phi_0$ and 1. The coefficient $s_0$ is the constrained specific storage coefficient.  The parameters $\lambda$ and $\mu$ are the Lam\'e parameters, and $K$ is the symmetric and positive-definite permeability tensor of the porous medium.  The system is forced by $\boldsymbol{f}$ and $\boldsymbol{g}$.

There are many alternative formulations of this system, and spatial discretizations are a much-studied topic~\cite{BOTH2022114183, kraus2024analysis}.
We follow the three-field model and its discretization studied in~\cite{kraus2024analysis}.
There, they discretize the displacement, $\boldsymbol{u}$, with an $H(\mathrm{div})$-conforming space, using stabilization terms to weakly enforce $H^1$ conformity, as in $H(\mathrm{div})-L^2$ discretizations of the Stokes equations~\cite{john2017divergence}.  The fluid variables, $\boldsymbol{w}$ and $\boldsymbol{p}$, are discretized with a stable $H(\mathrm{div}) \times L^2$ conforming pair.
Like in~\cite{kraus2024analysis}, we take the $H(\mathrm{div})$ space to be the Brezzi-Douglas-Marini space of some order $\ell$ over a triangulation of the domain $\Omega$, and the $L^2$-conforming space to be discontinuous polynomials of degree $\ell-1$ over the same triangulation.

In~\cite{kraus2024analysis}, the system is integrated in time by rewriting it as a first-order system in the standard way and using a continuous Petrov-Galerkin time stepping method, leading to parameter robustness and optimal order error estimates.
Here, we retain their spatial discretization and integrate in time with a fully implicit \RKN{} method.  Up to choices of quadrature, the Gauss-Legendre method is equivalent to their Galerkin-in-time method, but we obtained better results
using the RadauIIA method.
It is also possible to render the system as its first-order in time counterpart, integrating with a \RK{} method.

While the system in~\eqref{eq:biot_system} contains one second-order time derivative in two places, $\boldsymbol{u}_{tt}$, only first-order time derivatives of $\boldsymbol{w}$ and $p$ appear.  Just as in the case of DAEs, we can make use of the \RKN{} formalism even for terms in the equation that only have first-order time derivatives, simply by generalizing the replacement in~\eqref{eq:rknstageeq}, replacing all second-order time derivatives by the appropriate stage term, while replacing first-order time derivatives and zeroth-order time derivatives by the appropriate stage approximations as in~\eqref{eq:rknupdate}.  We note that this requires non-physical initial values for $\boldsymbol{w}_t$ and $p_t$, which we take to be zero.  While these values do enter into the stage approximations for $\boldsymbol{w}$ and $p$, they do so in a consistent manner, so that the values specified do not impact the resulting accuracy.

\section{Algebraic systems}
\label{sec:alg}
As demonstrated, implicit \RKN{} methods lead to large, stage-coupled algebraic systems, presenting greater practical difficulties than faced with single-stage or multi-step methods.
Recent years have seen major improvements in preconditioned iterative methods for these linear systems.
To help explicate these approaches, we 
restrict ourselves to the case of linear, constant-coefficient operators like the wave equation, but our implementation works much more generally.  Here, the linear system~\eqref{eq:blocksys} can be expressed using Kronecker products
\begin{equation}\label{eq:blocksys_kronecker}
  B \bfkappa \equiv \left( I \otimes M + \left( \Delta t \right)^2 \overline{A} \otimes K \right) \bfkappa = \bff.
\end{equation}

Preconditioners may largely be classified as to whether they segregate stages or treat them in a coupled fashion.
For first-order systems, we can point to early work~\cite{nilssen2011order, staff2006preconditioning} that employs the block diagonal or triangular part of the coupled system as a preconditioner.
This leads to a system requiring the inversion of diagonal blocks that allow reuse of a suitable strategy for single-stage methods.
However, these techniques lead to outer iteration counts that grow with stages.
\citet{masud2021new} introduce an alternative approach, where they first introduce a triangular approximation to the Butcher matrix based on the LDU factorization and thereby obtain a block triangular preconditioning matrix.
One can also obtain lower-triangular matrices that minimize the condition number $\kappa(L^{-1} A)$~\cite{staff2006preconditioning} or the norm $\| L^{-1} A - I \|_2$~\cite{rani2025efficient}.
This approach was extended to \RKN{} methods without first-order time derivatives in~\cite{clines2022efficient} and to \RK{} discretizations of first-order formulations of the wave equation in~\cite{rani2025efficient}.

In this case, any approximation $\overline{A} \approx \tilde{A}$, introduces an associated preconditioning matrix
\begin{equation}
  P = I \otimes M + \left( \Delta t \right)^2 \tilde{A} \otimes K,
\end{equation}
and typical concerns apply.  That is, $P$ must approximate $B$ well enough for the preconditioned iterative method to converge quickly, and $P$ must be significantly cheaper to invert onto a given vector than $B$.
For example, one obtains block-structured preconditioners by taking $\tilde{A}$ to be the diagonal or lower triangular part of $A$.  In the two-stage case, a lower-triangular approximation $\tilde{A}$ gives the preconditioning matrix
\begin{equation}
  P = \begin{bmatrix}
    M + \tilde{a}_{11} \left( \Delta t \right)^2 K
    & 0 \\
    \tilde{a}_{21} \left( \Delta t \right)^2 K &
    M + \tilde{a}_{22} \left( \Delta t \right)^2 K
  \end{bmatrix},
\end{equation}
and with an effective solver for the diagonal blocks (which are just single-stage discrete wave operators), this is readily inverted.  The approach of~\cite{clines2022efficient} takes the standard decomposition of $A = LDU$ and approximates $\tilde{A} = LD$, leading to helpful block-triangular structure, but also outer iteration counts that do not increase significantly with the number of \RKN{} stages.

On the other hand, \emph{monolithic multigrid} techniques use specialized relaxation schemes that couple all the stages together, applying multigrid principles directly to~\eqref{eq:blocksys_kronecker}.  Early work in~\cite{vanlent2005} showed that stage-coupled relaxation schemes could be very powerful, and this approach has been greatly extended and applied to complex fluids problems in~\cite{abu2022monolithic, farrell2021irksome}.
Suppose the semidiscrete method is posed on a finite element space $V_h$, and that some additive Schwarz decomposition of the form
\begin{equation}
  V_h = \sum_{i=1}^{N_s} V_i,
\end{equation}
where $N_s$ denotes the number of Schwarz subdomains,
provides an acceptable relaxation scheme for a single stage method.
The stage-coupled algebraic system is posed on the space
\begin{equation}
  \mathbf{V}_h^s = \prod_{j=1}^s V_h,
\end{equation}
and the monolithic multigrid algorithm again uses additive Schwarz relaxation, now based on the decomposition
\begin{equation}
  \mathbf{V}_h^s = \sum_{i=1}^{N_s} \mathbf{V}_i^s,
\end{equation}
where $\mathbf{V}_i^s = \prod_{j=1}^s V_i$.
A theory relating the convergence of multigrid methods based on monolithic relaxation to the convergence of equivalent multigrid for single-stage methods is given in~\cite{mmg}.

\section{Implementing \RKN{} methods in \Irksome}
\label{sec:imp}
Applying multi-stage \RKN{} methods to PDEs requires first writing out the stage-coupled variational problem to be solved at each time step.
Consider~\eqref{eq:2stagerknwave}, which could be written in UFL in a manner as shown in Figure~\ref{fig:2stagewaveufl}.
\begin{figure}
  \begin{lstlisting}
    Vbig = V * V
    k = Function(Vbig)
    k1, k2 = split(k)
    v1, v2 = TestFunctions(Vbig)
    u1 = u + c1 * dt * ut + dt**2 * (a11 * k1 + a12 * k2)
    u2 = u + c2 * dt * ut + dt**2 * (a21 * k1 + a22 * k2)
    F = (inner(k1, v1) * dx + inner(k2, v2) * dx
        + inner(grad(u1), grad(v1)) * dx + inner(grad(u2), grad(v2)) * dx)
\end{lstlisting}
\caption{UFL representation of the variational problem for a 2-stage implicit \RKN{} method.}
\label{fig:2stagewaveufl}
\Description{Code listing showing the UFL variational problem for \RKN{} stages.}
\end{figure}
Increasing the number of stages or the complexity of the variational problem, this system becomes tedious to express manually.  However, it can be programmatically obtained from the original problem and, so, is a highly suitable candidate for mechanization.
\Irksome{} now allows users to write a second-order semi-discrete variational problem in UFL, as shown in Figure~\ref{fig:fdnrkwave} in the introduction.
Internally, it processes these forms into the relevant systems for the \RKN{} stages, like the one shown in Figure~\ref{fig:2stagewaveufl}, and allows users access to a broad suite of solvers and preconditioners available at the PETSc/Firedrake interface~\cite{kirby2018solver}. 

\subsection{UFL Manipulation}

Hence, rather than implementing \RKN{} methods for equations of the form in~\eqref{eq:ode}, we consider variational evolution equations, seeking functions $y$ mapping the time interval $(0, T]$ into some finite-dimensional function space, $V_h$, such that  
\begin{equation}
  \label{eq:varode}
  \left( y_{tt}(t), v \right) = F(t, y, y_t; v)
\end{equation}
for all $v \in V_h$ where, as above, $\left( \cdot, \cdot \right)$ is some inner product on the space.  Even more generally, \RKN{} methods may be formulated for the fully nonlinear formulation
\begin{equation}
  \label{eq:genode}
  \mathcal{G}(t, y, y_t, y_{tt}; v) = 0.
\end{equation}
Although assertions about solvability and error estimates would surely require further assumptions, starting from this general form allows us to formulate \RKN{} methods for a very wide class of problems, including those examples discussed above.

In~\cite{farrell2021irksome}, we introduced an approach to UFL manipulation to obtain stage-coupled equations for standard \RK{} methods applied to systems of first-order ordinary differential equations.
As part of this work, we have greatly simplified our manipulation process and provided a clean extension to problems of second order.
In a UFL representation of~\eqref{eq:genode}, the unknown $y$ will be represented by some member of a \lstinline{FunctionSpace}, called \lstinline{V}.
This space itself may be a product of an arbitrary number of \lstinline{FunctionSpace} instances, each of which may be \lstinline{VectorFunctionSpace} or \lstinline{TensorFunctionSpace}.  The variable $y_t$ will be represented by a separate member of the same \lstinline{FunctionSpace}.
The variational problem for the stages is posed on a new \lstinline{FunctionSpace} called \lstinline{Vbig}, which is a \lstinline{MixedFunctionSpace} comprising the $s$-way product of the original \lstinline{V} with itself.

Then, the stage-coupled variational problem is obtained by summing over stages, making appropriate substitutions at each stage.
In particular, at stage $i$,
\begin{itemize}
\item The variable for $t$ is replaced with $t^n + c_i \Delta t$
\item The variable for $y$ is replaced by
  $y + c_i (\Delta t) y_t + (\Delta t)^2 \sum_{j=1}^s \overline{a}_{ij} \kappa_j$.
\item The variable for $y_t$ is replaced by
  $y_t + \Delta t \sum_{j=1}^s a_{ij} \kappa_j$
\item The variable for $y_{tt}$ is replaced by $\kappa_i$
\item The test function $v$ is replaced by the test function $v_i$, a particular piece of the \lstinline{TestFunction} associated with \lstinline{Vbig}.
\end{itemize}

To efficiently deal with the replacements on complicated expressions,
we updated the preprocessing of the semidiscrete formulation with a \lstinline{ufl.DAGTraverser}
to match the mechanism used on the spatial gradient and Gateaux derivative.
The \lstinline{ufl.DAGTraverser} expands the usual differentiation rules in order to ensure that
\lstinline{Dt} is only applied to the time-dependent coefficients of the problem.
For instance, \lstinline{Dt(div(y))} is preprocessed as \lstinline{div(Dt(y))}, and
\lstinline{Dt(sin(t))} is preprocessed as \lstinline{cos(t)}.
Hence, non-automonous forcing terms are symbolically differentiated.
To facilitate the preprocessing in the case where \lstinline{V} is compound (containing
mixed and/or vectorized function spaces), \lstinline{Dt} is commuted with indexing
and tensor operations.
For instance, \lstinline{Dt(as_vector([y[0], y[1]]))} becomes \lstinline{as_vector([Dt(y)[0], Dt(y)[1]])},
so that, in the end, we only require replacing the time derivative of the entire mixed function.

\subsection{Preconditioners}
Firedrake's rich interface to PETSc readily allows the implementation of the effective preconditioners discussed above in Section~\ref{sec:alg}.

Our implementation for the preconditioner proposed by~\citet{clines2022efficient} closely follows that for the analogous preconditioner for first order ODEs~\cite{kirby2024extending}.
Essentially, Irksome places information such as the semidiscrete UFL form and the tableau in a PETSc application context when the time stepping object is created.
Then, we implement an \lstinline{AuxiliaryOperatorPC}, as described in~\cite{kirby2018solver}, which is a harness for providing a preconditioning matrix derived from a different bilinear form than the system Jacobian.
In our case, the \lstinline{irksome.ClinesLD} preconditioner
extracts the semidiscrete UFL form and Nystr\"om tableau from an application context.
It computes the LDU factorization of $A$ and $\overline{A}$, reuses Irksome's UFL manipulation schemes to form the lower-triangular stage equations associated with this modified tableau.
The (block lower triangular) Jacobian of this form is then set as the auxiliary operator for the system.
It is important to note that the~\lstinline{AuxiliaryOperatorPC} does not specify how this matrix is to be (approximately) inverted, but simply exposes a new PETSc \lstinline{PC} object for the user to configure.
Figure~\ref{fig:clineslisting} shows a possible implementation.
These options use a matrix-free approach for the stage-coupled system matrix (only the action of the operator on a vector is computed).
Then, the \lstinline{ClinesLD} preconditioner is used, and it is configured to use a block triangular preconditioner via PETSc's \lstinline{fieldsplit}~\cite{brown2012composable}.
We then assemble the diagonal blocks of this preconditioner and apply one sweep of algebraic multigrid via \lstinline{hypre}~\cite{falgout2002hypre, VEHenson_UMYang_2002a} as an approximate inverse.

\begin{figure}
  \begin{lstlisting}
params = {
  "mat_type": "matfree",
  "snes_type": "ksponly",
  "ksp_type": "gmres",
  "ksp_rtol": 1.e-7,
  "pc_type": "python",
  "pc_python_type": "irksome.ClinesLD",
  "aux": {
    "ksp_type": "preonly",
    "pc_type": "fieldsplit",
    "pc_fieldsplit_type": "multiplicative",
    "fieldsplit": {
      "ksp_type": "preonly",
      "pc_type": "python",      
      "pc_python_type": "firedrake.AssembledPC",
      "aux_pc_type": "hypre"
    },
  }
}
\end{lstlisting}
  \caption{Typical options dictionary for implementing the Clines preconditioner for the wave equation}
  \label{fig:clineslisting}
  \Description{Code listing for the Clines preconditioner.}
\end{figure}

Monolithic geometric multigrid methods simply reuse existing Firedrake capabilities, and have been available since our first paper on Irksome~\cite{farrell2021irksome}.
Interfacing Firedrake with geometric multigrid in PETSc has long been available~\cite{mitchell2016high}.  Custom, patch-based relaxation schemes such as vertex star and Vanka relaxation are available in~\lstinline{PCPatch}~\cite{farrell2021pcpatch} and algebraic versions, such as~\lstinline{ASMStarPC} and~\lstinline{ASMVankaPC}.
No additional development was required to facilitate application of monolithic multigrid for \RKN{} methods.

\section{Numerical results}
\label{sec:numres}
Our numerical results are carried out on a 32-core AMD Ryzen Threadripper PRO 7975WX with 768GB of RAM.  The system is running Ubuntu 24.04. We use the configuration of Firedrake, \Irksome{}, PETSc, and other components of the software stack avalable through~\cite{zenodo/Zenodo-nirksome}.

\subsection{Wave equation}
For our first experiment, we consider the wave equation~\eqref{eq:semidwave}, posed on the unit cube $\Omega = [0, 1]^3$.
Our goal is to compare the effectiveness of implicit and explicit \RKN{} methods to a common explicit method, using central time differencing for $u_{tt}$.
Applying central time differencing to~\eqref{eq:waveODE} gives a simple problem to solve at each time step,
\begin{equation}
M \left( \frac{y^{n+1} - 2 y^n + y^{n-1}}{\Delta t^2} \right) = K y^n.
\end{equation}
Given the initial $y^0$ and a $y^{-1}$ (computed from $y^0$ and the initial data for the time derivative), one obtains each $y^{n+1}$ for $n\geq0$ from the two previous values by solving
\[
M y^{n+1} = M \left( 2 y^n - y^{n-1} \right) + \Delta t^2 K y^n.
\]
This requires some matrix-vector multiplication, vector addition, and the inversion of a mass matrix at each time step.
However, there are well-known techniques for obtaining diagonal approximations to $M$ so that no explicit inversion (via LU factorization or an iterative method) is required.

When $\Omega$ is decomposed into boxes, we use the standard $\mathcal{Q}_k$ space comprising tensor products of Lagrange polynomials associated with the Gauss-Lobatto points.
In this case, evaluating $M$ using tensor-product Gauss-Lobatto quadrature leads to a diagonal approximation to $M$.
With tetrahedral meshes and $k > 1$, there is no comparable set of points that parameterize the space $\mathcal{P}_k$ comprising polynomials of total degree $k$ while providing a suitably accurate quadrature rule.
However, the KMV elements~\cite{chin1999higher} enrich $\mathcal{P}_k$ with certain bubble functions and admit a set of nodal locations that coincide with a mass-lumping quadrature.  

We consider the following suite of experiments to measure the accuracy and efficiency of the two-stage Gauss-Legendre \RKN{} method and two explicit methods -- the centered difference and the classic Nystr\"om tableau~\eqref{eq:nystrom_scheme}.
We divide the unit cube into an $N \times N \times N$ mesh and apply a $\mathcal{Q}_k$ spatial discretization for $k=1,2$.
We begin with the initial condition of $u(x, y, z) = \sin(\pi x) \sin(\pi y) \sin(\pi z)$ and $u_t(x, y, z) = 0$ and integrate for two whole periods of oscillation, for a final time of $T_f=\tfrac{4}{\sqrt{3}}$.
With the implicit method, we use $\tfrac{N}{2}$ time steps per period, for $N$ total steps and a time step of $\Delta t = \tfrac{4}{N\sqrt{3}}$.
To solve the linear system, we use the stage-segregated Clines preconditioner, approximating the inverse of the diagonal blocks with one sweep of ML~\cite{MWGee_etal_2006a} using its default parameters.
For the central scheme, we compute an empirical bound for the time step based on a maximal eigenvalue estimated through a Krylov method.
For the explicit Nystr\"om method, we found empirically that taking $2kN$ time steps per period for $\Delta t = \frac{2}{kN\sqrt{3}}$ to be stable, but that a time step of half this was required to make the method have comparable accuracy to the other two.

\begin{figure}
  \begin{subfigure}[c]{0.49\textwidth}
    \caption{Steps taken}
    \label{fig:hexwavets}
    \begin{tikzpicture}[scale=0.6]
      \begin{axis}[xmode=log, ymode=log, log basis x={2}, ymax={3.e3},
          xlabel={$N$}, ylabel={$N_t$}, legend pos = south east]
        \addplot[black] table[x=N, y=nt, col sep = comma]{hex.rkn.deg3.cfl2.csv};
        \addlegendentry{Implicit};
        \addplot[red, mark=+]  table[x=N, y=nt, col sep = comma]{hex.central.deg1.csv};
        \addlegendentry{Central-$\mathcal{Q}_1$};
        \addplot[blue, mark=+] table[x=N, y=nt, col sep = comma]{hex.central.deg2.csv};
        \addlegendentry{Central-$\mathcal{Q}_2$};
        \addplot[red, dotted, mark=*]  table[x=N, y=nt, col sep = comma]{hex.nystrom.deg1.csv};
        \addlegendentry{Nystr\"om-$\mathcal{Q}_1$};
        \addplot[blue, dotted, mark=*] table[x=N, y=nt, col sep = comma]{hex.nystrom.deg2.csv};
        \addlegendentry{Nystr\"om-$\mathcal{Q}_2$};
      \end{axis}
    \end{tikzpicture}
  \end{subfigure}
  \begin{subfigure}[c]{0.49\textwidth}
    \caption{Accuracy}
    \label{fig:hexwaveacc}
    \begin{tikzpicture}[scale=0.6]
      \begin{axis}[xmode=log, ymode=log, log basis x={2},
          xlabel={$N$}, ylabel={$L^2$ error}]
        \addplot[dashed, red, mark=*] table[x=N, y=error, col sep = comma]{hex.rkn.deg1.cfl2.csv};
        \addplot[dashed, blue, mark=*] table[x=N, y=error, col sep = comma]{hex.rkn.deg2.cfl2.csv};
        \addplot[red, mark=+]  table[x=N, y=error, col sep = comma]{hex.central.deg1.csv};
        \addplot[blue, mark=+] table[x=N, y=error, col sep = comma]{hex.central.deg2.csv};
        \addplot[red, dotted, mark=*]  table[x=N, y=error, col sep = comma]{hex.nystrom.deg1.csv};
        \addplot[blue, mark=*, dotted] table[x=N, y=error, col sep = comma]{hex.nystrom.deg2.csv};          
      \end{axis}
    \end{tikzpicture}
  \end{subfigure}\\
  \begin{subfigure}[c]{0.49\textwidth}
    \caption{Runtime}
    \label{fig:hexwaveruntime}
    \begin{tikzpicture}[scale=0.6]
      \begin{axis}[xmode=log, ymode=log, log basis x={2}, xlabel={$N$},
          ylabel={Run time (s)},
          legend style={font=\LARGE}, legend pos= outer north east]
        \addplot[dashed, red, mark=*] table[x=N, y=time, col sep = comma]{hex.rkn.deg1.cfl2.csv};
        \addlegendentry{GL(2), $k=1$}          
        \addplot[dashed, blue, mark=*] table[x=N, y=time, col sep = comma]{hex.rkn.deg2.cfl2.csv};
        \addlegendentry{GL(2), $k=2$}
        \addplot[red, mark=+]  table[x=N, y=time, col sep = comma]{hex.central.deg1.csv};
        \addlegendentry{Central, $k=1$}
        \addplot[blue, mark=+] table[x=N, y=time, col sep = comma]{hex.central.deg2.csv};
        \addlegendentry{Central, $k=2$}
        \addplot[red, dotted, mark=*] table[x=N, y=time, col sep = comma]{hex.nystrom.deg1.csv};
        \addlegendentry{Nystr\"om, $k=1$}
        \addplot[blue, dotted, mark=*] table[x=N, y=time, col sep = comma]{hex.nystrom.deg2.csv};
        \addlegendentry{Nystr\"om, $k=2$}          
      \end{axis}
    \end{tikzpicture}
  \end{subfigure}
  \begin{subfigure}[c]{0.49\textwidth}
    \caption{GMRES performance for GL(2)}
    \label{fig:hexwaveits}
    \begin{tikzpicture}[scale=0.6]
      \begin{axis}[xmode=log, log basis x={2},
          xlabel={$N$}, ylabel={iterations}, legend pos = south east,
          ymin=0, ymax=15]
        \addplot[red, mark=+] table[x=N, y=its, col sep = comma]{hex.rkn.deg1.cfl2.csv};
        \addlegendentry{$\mathcal{Q}_1$};
        \addplot[blue, mark=+] table[x=N, y=its, col sep = comma]{hex.rkn.deg2.cfl2.csv};
        \addlegendentry{$\mathcal{Q}_2$};
      \end{axis}
    \end{tikzpicture}
  \end{subfigure}
  \caption{Comparing accuracy and run-time of two-stage Gauss-Legendre \RKN{} to explicit central differences and the Nystr\"om method for the wave equation discretized with $\mathcal{Q}_k$ elements on an $N \times N \times N$ hexahedral mesh.
    The wave equation is integrated over two whole periods.
    We report the time step size, accuracy at the final time, run-time, and (for iterative methods) the performance of our iterative solver.
  }
  \label{fig:wavehexperf}
  \Description{Plot showing the accuracy and efficiency of explicit and implicit solutions of the wave equation on hexahedra.}
\end{figure}
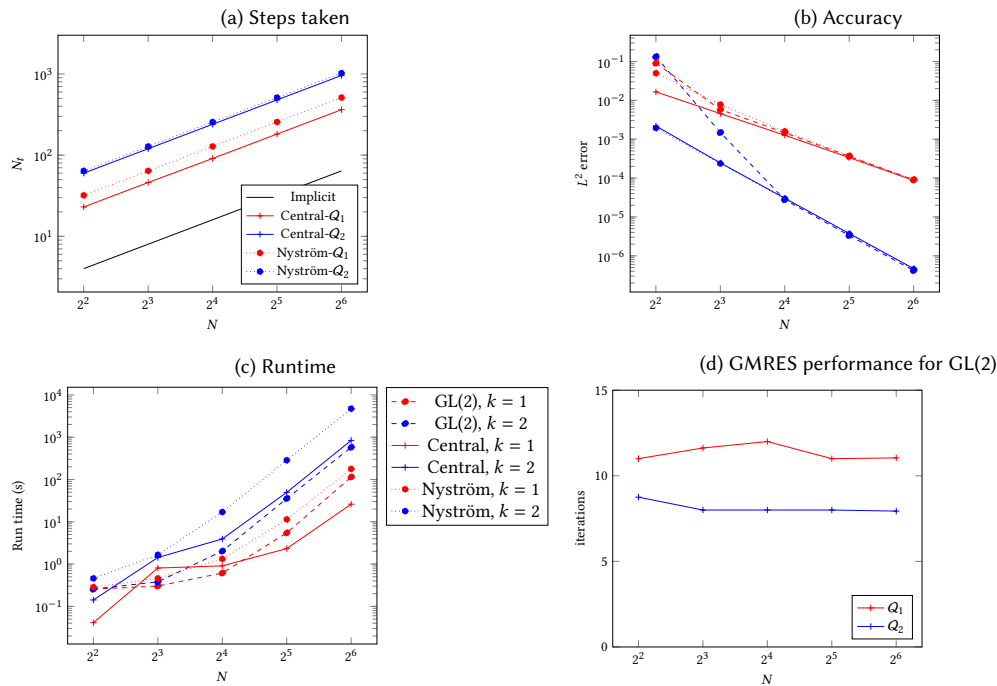

The various subplots in Figure~\ref{fig:wavehexperf} compare several aspects of the overall performance of our time stepping methods with linear and quadratic spatial discretization.
First, Figure~\ref{fig:hexwavets} plots the number of time steps taken for each our methods.
The implicit method (bottom black curve) takes many fewer time steps than either of the explicit schemes.
Then, we measure the accuracy of the method, plotting the $L^2$ norm of the error at the final time in Figure~\ref{fig:hexwaveacc}.
On coarse meshes, the GL(2) method produces some inaccuracy with $Q_2$, but on finer meshes this disappears.
Asymptotically, we see that all of our methods give about the same error, indicating that spatial error dominates temporal.
While this is expected for the formally fourth order GL(2) and Nystr\"om methods, it is perhaps surprising for the second-order central method.

Surprisingly, we see that for this particular problem, the implicit method performs quite well compared to the explicit ones.
It is faster than the Nystr\"om scheme for $\mathcal{Q}_1$ and is actually faster than both explicit methods with $\mathcal{Q}_2$.
A few factors contribute to this.  First, we are able to take very large time steps compared to the explicit methods.
Second, Figure~\ref{fig:hexwaveits} shows that the stage-segregated preconditioner is highly effective for this problem.
However, we note that we are able to reuse the assembled matrix and preconditioner at each time step.
For a nonlinear problem, the implicit method would incur the extra costs of a Newton iteration and updating the matrix and preconditioner, but we might still hope for large time steps and effective preconditioners.

\begin{figure}
  \begin{tikzpicture}[scale=0.6]
        \begin{axis}[xlabel={$t$}, ylabel={Energy}, legend pos = south east,
            ymin=1.82, ymax=1.88]
          \addplot[dashed] table[x=t, y=E, col sep = comma]{central_energy_1per_nx32_deg2.dat};
          \addlegendentry{Central};
          \addplot[mark=none] table[x=t, y=E, col sep = comma]{gl2_energy_1per_nx32_deg2.dat};
          \addlegendentry{GL(2)};
          \addplot[dotted] table[x=t, y=E, col sep = comma]{nystrom_energy_1per_nx32_deg2.dat};
          \addlegendentry{Nystr\"om};
        \end{axis}
  \end{tikzpicture}
    \caption{Discrete energy over one period of oscillation on a $32\times 32 \times 32$ hexahedral mesh using triquadratic elements. }
    \label{fig:energy}
      \Description{Plot showing energy conservation for explicit and implicit solutions of the wave equation.}
\end{figure}
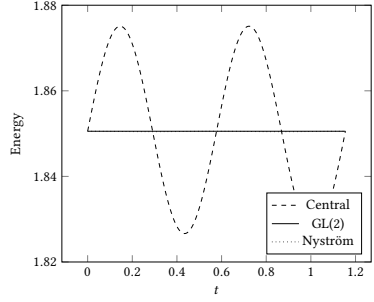

To illustrate the energy conservation of our time stepping methods, we integrate the wave equation for one period using quadratic elements on a $32\times 32 \times 32$ grid.
The energy as a function of time is plotted in Figure~\ref{fig:energy}.
The central difference method does not exactly conserve energy, but it oscillates around the initial value value.
In contrast, the implicit Gauss-Legendre method exactly (in practice, up to roundoff and solver tolerances) conserves the system energy exactly over time.
The Nystr\"om method does not have exact conservation, but only deviates slightly from the exact energy.

\subsection{Elastodynamic plates}
\begin{figure}[htbp]
  \begin{tikzpicture}[scale=0.5]
   \begin{loglogaxis}[xlabel={$N$},
        log basis x=2, legend style={font=\LARGE}, legend pos = outer north east
         ]
        \addplot table [x=N, y=errorL2, col sep=comma]{plate.nystrom.GaussLegendre-2.HCT-3.csv};   \addlegendentry{$\|u-u_h\|_{L^2}$}
        \addplot table [x=N, y=errorH1, col sep=comma]{plate.nystrom.GaussLegendre-2.HCT-3.csv};   \addlegendentry{$\|u-u_h\|_{H^1}$}
        \addplot table [x=N, y=errorH2, col sep=comma]{plate.nystrom.GaussLegendre-2.HCT-3.csv};   \addlegendentry{$\|u-u_h\|_{H^2}$}
        
        \addplot [domain=2^3:2^6] {2E-4/pow(x/2^6,2)} node[above, yshift=-2pt, midway, anchor=south west] {$h^{2}$};
        \addplot [domain=2^3:2^6] {4E-7/pow(x/2^6,3)} node[above, yshift=-2pt, midway, anchor=south west] {$h^{3}$};
        \addplot [domain=2^3:2^6] {1.5E-8/pow(x/2^6,4)} node[above, yshift=-2pt, midway, anchor=south west] {$h^{4}$};
   \end{loglogaxis}
  \end{tikzpicture}
        
\caption{Error norms computed at the final time for the elastodynamic plate discretized with cubic HCT macroelements and two-stage Gauss-Legendre.}
  \label{fig:platerr}
\Description{Figure showing accuracy for dynamic plate problem.}
\end{figure}
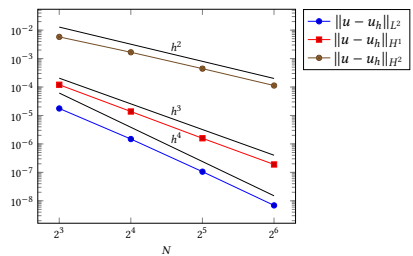

Next, we consider the elastodynamic plate problem~\eqref{eq:dynplate} posed on $\Omega = [0, 1]^2$, choosing the initial data so that the exact solution satisfies
\[
u(t, x, y) = \left[ x (1-x) y (1-y) \right]^2 \cos(\pi t).
\]
We divide $\Omega$ into an $N \times N$ mesh of squares, each subdivided into two right triangles, and then discretize our problem in with third-order Hsieh-Clough-Tocher elements~\cite{clough1965finite}, which were recently added to our code stack~\cite{brubeck2025fiat}.
Because of the energy conservation property~\eqref{eq:dynplateenergy}, we use a two-stage Gauss-Legendre scheme in time and integrate to $t=1$.
Figure~\ref{fig:platerr} shows the error versus mesh spacing, taking $\Delta t = \tfrac{1}{N}$.
We see nearly optimal convergence rates of second order in $H^2$, third order in $H^1$, and fourth order in $L^2$.

As a remark, obtaining optimal convergence rates in $L^2$ turned out to be somewhat delicate.
For the standard wave equation, \citet{dupont1973l2} showed that the $L^2$ error depends on the $H^1$ error in approximating the initial condition, requiring one to compute an elliptic rather than $L^2$ projection of the initial data.
After initially observing suboptimal convergence rates in $L^2$ for the dynamic wave equation, we replaced the $L^2$ projection of the initial condition with a biharmonic projection in an apparently analogous situation, which restored optimal convergence rates.

\begin{figure}
  \begin{subfigure}[c]{0.33\textwidth}
  \begin{tikzpicture}[scale=0.5]    
   \begin{loglogaxis}[xlabel={$N$}, ylabel={$\Delta t$},
        log basis x=2, legend style={font=\LARGE}, legend pos = north west
         ]
        \addplot table [x=N, y=timesteps, col sep=comma]{plate.central.HCT-3.csv};   \addlegendentry{Central}
        \addplot table [x=N, y=timesteps, col sep=comma]{plate.nystrom.GaussLegendre-2.HCT-3.csv};   \addlegendentry{Implicit}
   \end{loglogaxis}
  \end{tikzpicture}
  \caption{Number of time steps}
  \end{subfigure}
  \begin{subfigure}[c]{0.33\textwidth}
    \begin{tikzpicture}[scale=0.5]
      \begin{loglogaxis}[xlabel={$N$}, ylabel={Run time(s)},
          log basis x=2, legend style={font=\LARGE}, legend pos = north west
        ]
        \addplot table [x=N, y=runtime, col sep=comma]{plate.central.HCT-3.csv};   \addlegendentry{Central}
        \addplot table [x=N, y=runtime, col sep=comma]{plate.nystrom.GaussLegendre-2.HCT-3.csv};   \addlegendentry{RKN-GL(2)}
        \addplot table [x=N, y=runtime, col sep=comma]{plate.rk.GaussLegendre-2.HCT-3.csv};   \addlegendentry{RK-GL(2)}  
      \end{loglogaxis}
    \end{tikzpicture}
    \caption{Run-time}
  \end{subfigure}
  \begin{subfigure}[c]{0.33\textwidth}
    \begin{tikzpicture}[scale=0.5]
      \begin{semilogxaxis}[xlabel={$N$}, ylabel={Avg. Iterations}, ymin=0,
          log basis x=2, legend style={font=\LARGE}, legend pos = north east
        ]
        \addplot table [x=N, y expr=\thisrow{iterations}/\thisrow{timesteps}, col sep=comma]{plate.nystrom.GaussLegendre-2.HCT-3.csv};   \addlegendentry{RKN-GL(2)}
        \addplot table [x=N, y expr=\thisrow{iterations}/\thisrow{timesteps}, col sep=comma]{plate.rk.GaussLegendre-2.HCT-3.csv};   \addlegendentry{RK-GL(2)}  
      \end{semilogxaxis}
    \end{tikzpicture}
    \caption{GMRES Iterations}
  \end{subfigure}  
  \caption{Performance comparison of explicit central differences and \RKN{} and \RK{} methods for the dynamic plate problem discretized with HCT elements.
    We use an empirically-determined stable time step for the explicit method and $\Delta t = \tfrac{1}{N}$ for the implicit method.  Both \RKN{} formulation and first-order system \RK{} form use the two-stage Gauss-Legendre method.}
  \label{fig:dynplate}
  \Description{Figure comparing performance of explicit versus implicit methods for dynamic plate problem.}
\end{figure}
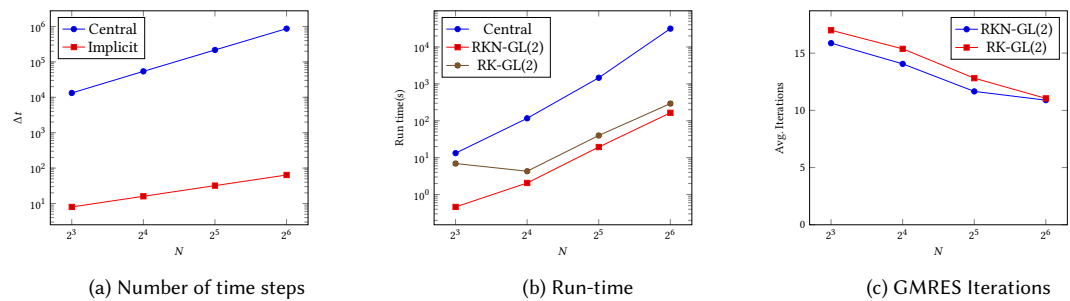

Although~\eqref{eq:dynplate} is a wave-type equation, the higher order spatial derivatives make the problem poorly-suited for explicit time stepping.
If we apply the central difference scheme considered above for the wave equation, the time step will be governed by the dominant generalized eigenvalue of
\begin{equation}
  K \mathbf{x} = \lambda M \mathbf{x},
\end{equation}
where $K$ is the biharmonic stiffness matrix and $M$ is the mass matrix associated with the HCT space.
On an $N \times N$ mesh, we estimated this dominant eigenvalue using a Krylov method and set the time step as $\Delta t = \tfrac{2}{\sqrt{\lambda_{max}}}$.
Figure~\ref{fig:dynplate} shows the the stable time step required for explicit time stepping as a function of mesh size and compares the run-time required to integrate to time $T=1$ using both explicit and implicit methods.
With $C^1$ elements like HCT, mass lumping is not possible, and we invert the mass matrix at each time step using conjugate gradients preconditioned with an incomplete Cholesky factorization.
In both cases, we set \lstinline{ksp_guess_type} to be \lstinline{fischer} to use the technique from~\cite{fischer1998projection}.
When solving a linear system repeatedly with different right-hand sides, this technique first approximates the solution in the span of prior solutions and then uses the Krylov method to correct the result.

Figure~\ref{fig:dynplate} also includes data comparing the performance of \RKN{} time stepping (denoted by RKN) and \RK{} time stepping on the equivalent first-order system (denoted by RK).  We see that there is a clear cost in run-time using the equivalent first-order system, of almost a factor of 2 for most values of $N$.  These results were computed using geometric multigrid W(2,2)-cycles to precondition FGMRES.  To achieve this performance we use vertex star patch-based relaxation coupling the \RKN{} (or \RK{}) stages, accelerated with GMRES.  The resulting iterations, averaged over the full run, are shown in the right panel of Figure~\ref{fig:dynplate}, showing that there is no loss in performance using \RKN{} methods.

\subsection{Poroelasticity}

For our final example, we consider the dynamic equations of Biot poroelasticity, following the formulation given in~\cite{kraus2024analysis}, described above.  We first consider results using the method of manufactured solutions, modifying the quasi-static test case from~\cite{fu2019-aa}.  Here, we consider the solution
\begin{align*}
\boldsymbol{u}(x,y,t) & = \sin(t)\left ( \begin{array}{c}\sin(\pi y) \left ( -\cos(\pi x) + \frac{1}{\mu + \lambda}\sin(\pi x) \right )\\
\sin(\pi x) \left ( \cos(\pi y) + \frac{1}{\mu + \lambda}\sin(\pi y) \right )\\\end{array} \right ),\\
p(x,y,t) & = \sin(t)\sin(\pi x)\sin(\pi y),\\
\boldsymbol{w}(x,y,t) &=-k\nabla p,
\end{align*}
computing consistent forcing terms in~\eqref{eq:biot_system}, and using these expressions to prescribe initial conditions and strong boundary conditions for the normal components of $\boldsymbol{u}$ and $\boldsymbol{w}$ on the boundary of a unit-square domain.  The tangential component of $\boldsymbol{u}$ is weakly enforced to match this prescription through the SIPG terms.  In these experiments, we will vary the Poisson ratio, $\nu$, while fixing the Young's modulus, $E = 3 \times 10^4$, yielding elasticity parameters
\[
\mu = \frac{E}{2 + 2\nu} \text{ and }
\lambda = \frac{E\nu}{(1- 2\nu)(1+\nu)}.
\]
We note that there is particular interest in computing solutions to~\eqref{eq:biot_system} in the incompressible limit as $\nu \rightarrow 0.5$, where the $\nabla(\nabla\cdot \boldsymbol{u})$ term is dominant, and the discretization considered here should be robust in this limit.  We fix the remaining parameters with $\alpha = 1$, $K = 10^6 I$, $s_0 = 10^{-6}$, $\rho_f = 1000$, $\rho_s = 500$, $\phi_0 = 0.1$, $\rho_w = \rho_f/\phi_0 = 10^4$, and $\overline{\rho} = (1-\phi_0)\rho_s + \phi_0\rho_w = 550$.

For the manufactured solution test, we consider spatial grids constructed by applying a fixed number of uniform refinement steps to a base mesh of $16\times 16$ quadrilaterals each cut into two triangles from the top left to bottom right corner. We integrate from time 0 to time 1 using a prescribed number of time steps, $n_t$. We vary both the order, $\ell$, of the mixed finite element space (using BDM($\ell$) elements for $\boldsymbol{u}$ and $\boldsymbol{w}$ and discontinuous $\mathcal{P}_{\ell-1}$ elements for $p$, and the number of stages in the \RKN{} method, considering methods using both the Gauss-Legendre and Radau IIA tableaux.  We solve the resulting linear systems at each time step using the monolithic multigrid scheme of~\cite{scott_biot}, which applies Vanka-style relaxation to the fully coupled discretization.

Convergence theory for both the quasi-static and dynamic cases is well-understood, as Biot poroelasticity forms a DAE system of order 1 (see~\cite{scott_biot} for details).  In this setting, it is easy to verify that the expected convergence for a stiffly accurate time stepper, such as Radau IIA, should be the formal order of the scheme for all components of the system.  For a non stiffly accurate scheme, such as Gauss-Legendre, we expect to achieve the formal order of the scheme only for the pressure term (which appears with a time derivative), and convergence at the stage order plus 1 for the displacement and velocity terms.  See~\cite[Section VI, Theorem 3.23]{wanner1996solving} for details on these rates.  We note that these convergence estimates do not directly apply to the \RKN{} schemes used here, nor do they account for the spatial discretization; see~\cite{scott_biot} for more complete theory.

\begin{figure}
  \begin{subfigure}[c]{0.33\textwidth}
  \begin{tikzpicture}[scale=0.5]    
   \begin{loglogaxis}[xlabel={$n_x$}, ylabel={Error},
        log basis x=2, legend style={font=\LARGE}, legend pos = south west
         ]
        \addplot table [x=Nx, y=product, col sep=comma]{biot_0.49_Nystrom_Radau.txt};   \addlegendentry{RKN-Radau(2)}
        \addplot table [x=Nx, y=product, col sep=comma]{biot_0.49_Nystrom_GL.txt};   \addlegendentry{RKN-GL(2)}
        \addplot table [x=Nx, y=product, col sep=comma]{biot_0.49_RK_GL.txt};   \addlegendentry{RK-GL(2)}
        \addplot [domain=2^5:2^7] {4/pow(x/2^5,1)} node[above, yshift=-2pt, midway, anchor=south west] {$h^{1}$};
        \addplot [domain=2^5:2^7] {1.5/pow(x/2^5,2)} node[above, yshift=-2pt, midway, anchor=south west] {$h^{2}$};
   \end{loglogaxis}
  \end{tikzpicture}
  \caption{Product-norm error}\label{fig:biot_mms_049_error}
  \end{subfigure}
  \begin{subfigure}[c]{0.33\textwidth}
    \begin{tikzpicture}[scale=0.5]
      \begin{loglogaxis}[xlabel={$n_x$}, ylabel={Run time(s)},
          log basis x=2, legend style={font=\LARGE}, legend pos = north west
        ]
        \addplot table [x=Nx, y=time, col sep=comma]{biot_0.49_Nystrom_Radau.txt};   \addlegendentry{RKN-Radau(2)}
        \addplot table [x=Nx, y=time, col sep=comma]{biot_0.49_Nystrom_GL.txt};   \addlegendentry{RKN-GL(2)}
        \addplot table [x=Nx, y=time, col sep=comma]{biot_0.49_RK_GL.txt};   \addlegendentry{RK-GL(2)}
        \addplot [domain=2^5:2^7] {20*pow(x/2^5,3)} node[below, yshift=-2pt, left, anchor=south west] {$n_x^3$};
      \end{loglogaxis}
    \end{tikzpicture}
    \caption{Run-time}\label{fig:biot_mms_049_time}
  \end{subfigure}
  \begin{subfigure}[c]{0.33\textwidth}
    \begin{tikzpicture}[scale=0.5]
      \begin{semilogxaxis}[xlabel={$n_x$}, ylabel={Avg. Iterations}, ymin=0,
          log basis x=2, legend style={font=\LARGE}, legend pos = south east
        ]
        \addplot table [x=Nx, y=its, col sep=comma]{biot_0.49_Nystrom_Radau.txt};   \addlegendentry{RKN-Radau(2)}
        \addplot table [x=Nx, y=its, col sep=comma]{biot_0.49_Nystrom_GL.txt};   \addlegendentry{RKN-GL(2)}
        \addplot table [x=Nx, y=its, col sep=comma]{biot_0.49_RK_GL.txt};   \addlegendentry{RK-GL(2)}
      \end{semilogxaxis}
    \end{tikzpicture}
    \caption{Preconditioned GMRES Iterations}\label{fig:biot_mms_049_its}
  \end{subfigure}  
  \caption{Performance comparison of \RKN{} and \RK{} methods for the Biot manufactured solution problem.  We fix the discretization order to be $\ell = 2$ and use 2 stages for the time steppers, considering simulations with an $n_x \times n_y$ finest spatial grid for $n_x = n_y$ with $h_t = 4/n_x$.  Here, we fix the value of $\nu = 0.49$.}
  \label{fig:biot_mms_049}
  \Description{Figure comparing performance of \RKN{} and \RK{} time stepping for Biot problem.}
\end{figure}
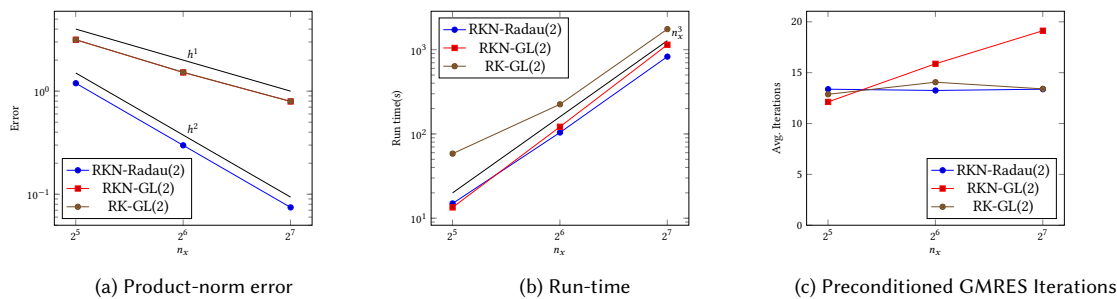

\Cref{fig:biot_mms_049} presents results for $\nu = 0.49$ with $\ell = 2$ and 2-stage time steppers.  Here, we consider simultaneous refinement of the spatial and temporal meshes, solving the problem on $n_x \times n_y$ meshes with $n_x = n_y$ and $n_t = n_x/4$.  At left, we measure the error in the product norm considered in~\cite{JAdler_etal_2021a}, given by
\[
\|(\boldsymbol{u},\boldsymbol{w},p)\|_\ast^2 = 2\mu \|\epsilon(\boldsymbol{u})\|^2 + \lambda \|\nabla\cdot\boldsymbol{u}\|^2 + 10^{-6}h_t\|\boldsymbol{w}\|^2 + c_ph_t^2\|\nabla\cdot\boldsymbol{w}\|^2 + c_p^{-1}\|p\|^2,
\]
where $c_p = \left(\frac{1}{\mu + \lambda} + s_0 \right)^{-1}$.  While this product norm (used to analyse a low-order Taylor-Hood spatial discretization of the quasistatic equations with implicit Euler time stepping) is not a perfect fit for our current setting, it does a good job of demonstrating the convergence of the discretization.  We note that since we use discontinuous piecewise linear functions for $p$ (which appears only with an $L^2$ norm in the product norm) and quadratic BDM elements for $\boldsymbol{u}$ and $\boldsymbol{w}$ (which appear with derivatives in the product norm), this limits our expected spatial convergence to second order.  In~\cref{fig:biot_mms_049_error}, we see that the \RKN{} method using the RadauIIA Butcher tableau converges with second order as expected, while the \RK{} and \RKN{} methods using Gauss time stepping converges only at first order.  In~\cref{fig:biot_mms_049_time}, we see that run-time for all three approaches scale roughly like $n_x^3$, indicating roughly linear scaling with total problem size ($\mathcal{O}(n_x^2)$) times the number of time steps (which scale like $\mathcal{O}(n_x)$), but that there is a clear savings in the use of \RKN{} approaches over direct \RK{} time steppers.  Finally, \cref{fig:biot_mms_049_its} shows the average number of preconditioned GMRES iterations per time step for the three methods.  Here, we see some growth with \RKN{} using the Gauss Butcher tableau, accounting for its slightly increasing cost in the middle graph, but roughly stable iteration counts for the other two approaches.

\begin{figure}
  \begin{subfigure}[c]{0.33\textwidth}
  \begin{tikzpicture}[scale=0.5]    
   \begin{loglogaxis}[xlabel={$n_x$}, ylabel={Error},
        log basis x=2
         ]
        \addplot table [x=Nx, y=product, col sep=comma]{biot_0.49_Nystrom_Radau_1stage.txt};
        \addplot table [x=Nx, y=product, col sep=comma]{biot_0.49_Nystrom_Radau_2stages.txt};
        \addplot table [x=Nx, y=product, col sep=comma]{biot_0.4999_Nystrom_Radau_1stage.txt};
        \addplot table [x=Nx, y=product, col sep=comma]{biot_0.4999_Nystrom_Radau_2stages.txt};
        \addplot [domain=2^5:2^7] {40/pow(x/2^5,1)} node[above, yshift=-2pt, midway, anchor=south west] {$h^{1}$};
        \addplot [domain=2^5:2^7] {0.5/pow(x/2^5,3)} node[above, yshift=-2pt, midway, anchor=south west] {$h^{3}$};
   \end{loglogaxis}
  \end{tikzpicture}
  \caption{Product-norm error}\label{fig:biot_mms_04999_error}
  \end{subfigure}
  \begin{subfigure}[c]{0.33\textwidth}
    \begin{tikzpicture}[scale=0.5]
      \begin{loglogaxis}[xlabel={$n_x$}, ylabel={Run time(s)},
          log basis x=2
        ]
        \addplot table [x=Nx, y=time, col sep=comma]{biot_0.49_Nystrom_Radau_1stage.txt};
        \addplot table [x=Nx, y=time, col sep=comma]{biot_0.49_Nystrom_Radau_2stages.txt};
        \addplot table [x=Nx, y=time, col sep=comma]{biot_0.4999_Nystrom_Radau_1stage.txt};
        \addplot table [x=Nx, y=time, col sep=comma]{biot_0.4999_Nystrom_Radau_2stages.txt};
        \addplot [domain=2^5:2^7] {10*pow(x/2^5,3)} node[below, yshift=-2pt, left, anchor=south west] {$n_x^3$};
      \end{loglogaxis}
    \end{tikzpicture}
    \caption{Run-time}\label{fig:biot_mms_04999_time}
  \end{subfigure}
  \begin{subfigure}[c]{0.33\textwidth}
    \begin{tikzpicture}[scale=0.5]
      \begin{semilogxaxis}[xlabel={$n_x$}, ylabel={Avg. Iterations}, ymin=0,
          log basis x=2, legend style={font=\LARGE}, legend pos = south east
        ]
        \addplot table [x=Nx, y=its, col sep=comma]{biot_0.49_Nystrom_Radau_1stage.txt};   \addlegendentry{$\nu = 0.49$,Radau(1)}
        \addplot table [x=Nx, y=its, col sep=comma]{biot_0.49_Nystrom_Radau_2stages.txt};   \addlegendentry{$\nu = 0.49$,Radau(2)}
        \addplot table [x=Nx, y=its, col sep=comma]{biot_0.4999_Nystrom_Radau_1stage.txt};   \addlegendentry{$\nu = 0.4999$,Radau(1)}
        \addplot table [x=Nx, y=its, col sep=comma]{biot_0.4999_Nystrom_Radau_2stages.txt};   \addlegendentry{$\nu = 0.4999$,Radau(2)}
      \end{semilogxaxis}
    \end{tikzpicture}
    \caption{Preconditioned GMRES Iterations}\label{fig:biot_mms_04999_its}
  \end{subfigure}  
  \caption{Performance of \RKN{} method for the Biot manufactured solution problem with Radau time steppers and varying $\nu$.  We consider two orders of discretization, using 1 stage and $\ell = 1$ (denoted Radau(1)), or 2 stages and $\ell = 3$ (denoted Radau(3)), considering simulations with an $n_x \times n_y$ finest spatial grid for $n_x = n_y$ with $h_t = 4/n_x$.}
  \label{fig:biot_mms_04999}
  \Description{Figure comparing performance of \RKN{} time stepping for Biot problem in incompressible limit.}
\end{figure}
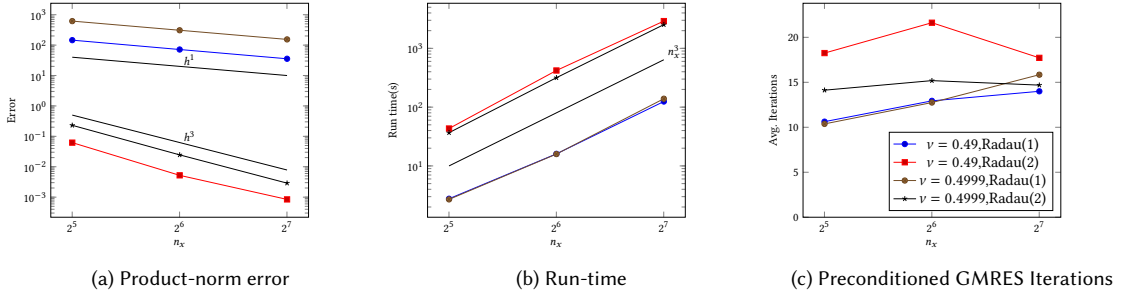

Next, we consider the \RKN{} method using Radau IIA time stepping for a more incompressible case, with $\nu = 0.4999$, as we simultaneously increase both the polynomial degree and the number of stages in the \RKN{} method.  \Cref{fig:biot_mms_04999} presents the same metrics as above for both $\nu = 0.49$ and $\nu = 0.4999$.  Here, we consider two combinations of discretization order, informed by the DAE analysis.  As a low-order discretization, we use 1 stage for time stepping (implicit Euler) in combination with $\ell = 1$.  Here, we expect first-order convergence, as we see for both values of $\nu$ in~\cref{fig:biot_mms_04999_error}.  As a higher-order combination, we consider the 2 stage Radau method, in combination with $\ell = 3$.  Here, we expect third-order convergence from both spatial and temporal discretization, which we also see in~\cref{fig:biot_mms_04999_error}.  While there is a visible degradation in errors between the more compressible ($\nu = 0.49$) and more incompressible ($\nu = 0.4999$) cases, this does not appear to depend on grid size, and is less than an order of magnitude in all cases.  Runtime differences between the two cases are, in contrast, negligible, as shown in~\cref{fig:biot_mms_04999_time}.  In~\cref{fig:biot_mms_04999_its}, we see no systematic differences between the preconditioned GMRES iteration counts between all cases.

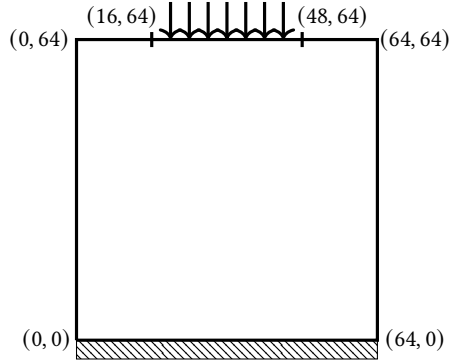
\begin{figure}
        \centering
        \begin{tikzpicture}[scale=1]

          \draw[very thick] (0,0) -- (4,0) -- (4, 4) -- (0, 4) -- (0,0);
          \draw[pattern=north west lines] (0,0) rectangle (4,-0.25);
          \draw[very thick] (1,3.9) -- (1, 4.1);
          \draw[very thick] (3,3.9) -- (3, 4.1);
          \foreach \x in {1.25, 1.5, ..., 2.75} {\draw[very thick, ->] (\x, 4.5) -- (\x, 4);};
          \node[align=left] at (-0.4,0) {$(0,0)$};
          \node[align=right] at (4.5,0) {$(64,0)$};
          \node[align=right] at (4.5,4) {$(64,64)$};
          \node[align=left] at (-0.5,4) {$(0,64)$};
          \node[align=left] at (0.6,4.3) {$(16,64)$};
          \node[align=right] at (3.4,4.3) {$(48,64)$};
        \end{tikzpicture}
        \caption{Schematic for Biot footing problem.  A downward traction is posed on the segmented below the arrows.}
        \Description{Picture of the domain for Biot footing problem.}
        \label{fig:footing}
\end{figure}

The second poroelasticity problem that we consider is a footing problem, adapted from~\cite{BOTH2022114183}, considering the square domain, $[0,64]^2$, with boundary conditions as pictured in~\cref{fig:footing}.  Zero flow boundary conditions are imposed on $\boldsymbol{w}$ on the middle segment of the top face, with zero pressure boundary conditions imposed on the complement of that segment.  Zero displacement boundary conditions are enforced on $\boldsymbol{u}$ on the bottom face (strongly for the normal component, weakly for the tangential component), while natural boundary conditions are imposed on the displacement on all other faces.  On the middle segment of the top face, a traction force in the downward direction of magnitude $10^5 t$ is imposed, with zero traction force on all other faces with natural displacement boundary conditions.  We again consider integration for $0 \leq t \leq 1$ using $n_t$ time steps, on grids that are uniform refinements of a $16\times 16$ base mesh constructed as above (but on the $[0,64]^2$ domain).  We again use the monolithic multigrid method from~\cite{scott_biot} to solve the resulting linear systems.

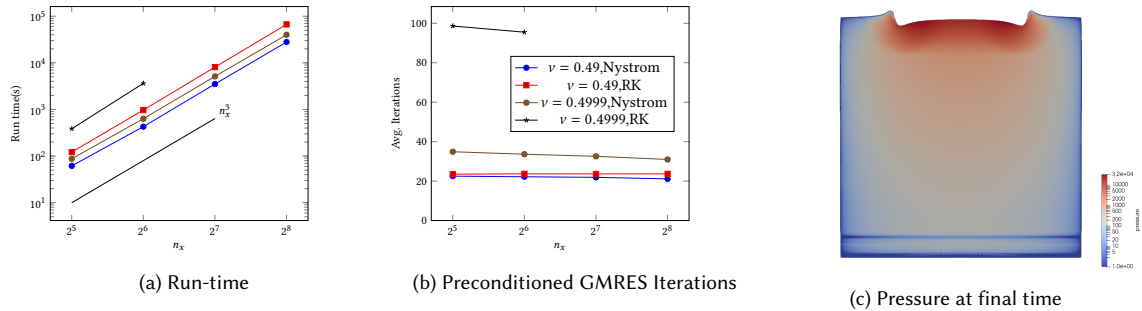
\begin{figure}
  \begin{subfigure}[c]{0.33\textwidth}
    \begin{tikzpicture}[scale=0.5]
      \begin{loglogaxis}[xlabel={$n_x$}, ylabel={Run time(s)},
          log basis x=2
        ]
        \addplot table [x=Nx, y=time, col sep=comma]{biot_0.49_footing_nystrom.txt};
        \addplot table [x=Nx, y=time, col sep=comma]{biot_0.49_footing_rk.txt};
        \addplot table [x=Nx, y=time, col sep=comma]{biot_0.4999_footing_nystrom.txt};
        \addplot table [x=Nx, y=time, col sep=comma]{biot_0.4999_footing_rk.txt};
        \addplot [domain=2^5:2^7] {10*pow(x/2^5,3)} node[below, yshift=-2pt, left, anchor=south west] {$n_x^3$};
      \end{loglogaxis}
    \end{tikzpicture}
    \caption{Run-time}\label{fig:biot_footing_time}
  \end{subfigure}
  \begin{subfigure}[c]{0.33\textwidth}
    \begin{tikzpicture}[scale=0.5]
      \begin{semilogxaxis}[xlabel={$n_x$}, ylabel={Avg. Iterations}, ymin=0,
          log basis x=2, legend style={font=\LARGE,at={(axis cs:56,45)},anchor=south west},
        ]
        \addplot table [x=Nx, y=its, col sep=comma]{biot_0.49_footing_nystrom.txt};   \addlegendentry{$\nu = 0.49$,Nystrom}
        \addplot table [x=Nx, y=its, col sep=comma]{biot_0.49_footing_rk.txt};   \addlegendentry{$\nu = 0.49$,RK}
        \addplot table [x=Nx, y=its, col sep=comma]{biot_0.4999_footing_nystrom.txt};   \addlegendentry{$\nu = 0.4999$,Nystrom}
        \addplot table [x=Nx, y=its, col sep=comma]{biot_0.4999_footing_rk.txt};   \addlegendentry{$\nu = 0.4999$,RK}
      \end{semilogxaxis}
    \end{tikzpicture}
    \caption{Preconditioned GMRES Iterations}\label{fig:biot_footing_its}
  \end{subfigure}  
  \begin{subfigure}[c]{0.33\textwidth}
    \includegraphics[width=\textwidth]{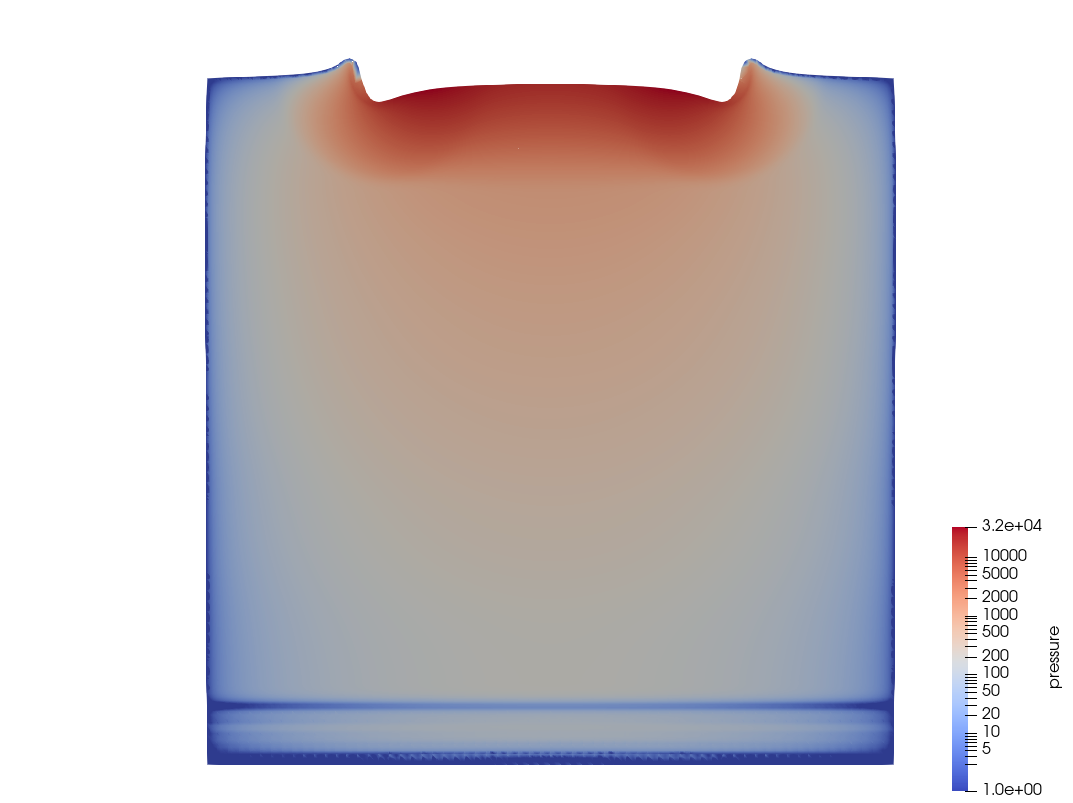}
    \caption{Pressure at final time}\label{fig:biot_footing_under_pressure}
  \end{subfigure}  
  \caption{Performance of \RKN{} and \RK{} methods for the Biot footing problem with 2-stage Radau time steppers and $\ell = 3$ and varying $\nu$.  Simulations use an $n_x \times n_y$ finest spatial grid for $n_x = n_y$ with $h_t = 4/n_x$.}
  \label{fig:biot_footing}
  \Description{Figure comparing performance of \RKN{} and \RK{} time stepping for Biot footing problem.}
\end{figure}

\Cref{fig:biot_footing} details the performance of the \RKN{} and \RK{} time steppers for solving this problem using a 2-stage RadauIIA time stepper and $\ell=3$, for $\nu = 0.49$ and $0.4999$, using FGMRES(50) as the outer Krylov solver.  In~\cref{fig:biot_footing_time}, we see that the CPU time to solution scales cubically for both methods applied to both problems, with the \RKN{} time stepping requiring less time than \RK{} methods.  For $\nu = 0.49$, the \RKN{} method is consistently about two times faster than \RK{}, increasing to abou 2.4 times faster for $n_x = 256$.  For $\nu = 0.4999$, while the \RKN{} method is slightly slower than it is for $\nu = 0.49$, it is substantially faster than \RK{}, for which the linear solver fails to converge beyond $n_x = 64$, even allowing up to a total of 1000 iterations of FGMRES(50).  Even for the cases where the \RK{} discretization does yield a solution, it is substantially less efficient, requiring over 5 times the computational time for $n_x = 64$.  This degradation in performance is partly explained in~\cref{fig:biot_footing_its}, which shows that the linear solver performs quite poorly for the \RK{} discretization with $\nu = 0.4999$, while it offers similar performance for all the other 3 problems and methods shown.  In particular, we see very little difference in iteration counts for the two discretizations at $\nu = 0.49$ (indicating that the additional time for the \RK{} discretization comes almost entirely from the larger system of equations), and that the variation in \RKN{} timing with $\nu$ is due to an increase in iterations in the nearly incompressible case.  \Cref{fig:biot_footing_under_pressure} shows the fluid pressure at the final time of $t=1$, noting similar results as those shown in~\cite{BOTH2022114183}.

\section{Conclusions}
\label{sec:conc}
The classical but rarely-implemented \RKN{} methods offer effective tools for directly discretizing PDEs with second-order time derivatives.
Through Irksome and Firedrake, we can automate the construction of the rather complex stage-coupled variational problems to be solved at each time step, and have deployed a range of effective solvers for implicit problems.
Appropriately tuned, implicit methods can be competitive for wave equations traditionally integrated with explicit schemes.  Applying \RKN{} directly to second-order problems rather than applying the underlying \RK{} schemes to equivalent first-order systems offers nontrivial performance advantages. 

Wave equations and other problems with second-order time derivatives arise in many common applications, such as geophysical imaging and structural mechanics.  Possible avenues for future work include using \Irksome{} to model these systems in scientific and industrial applications, and (in combination with Firedrake's ability to compute adjoints) to solve the corresponding inverse problems.  We also note significant possibilities for further development of preconditioning techniques for both first- and second-order problems that is facilitated by having the tool presented here to automatically generate \RKN{} discretizations for second-order PDEs.

\begin{acks}
  The work of P.~D.~B.~was supported by EPSRC grant EP/W026260/1.
  The work of R.~C.~K.~was partially supported by NSF 2410408.
  The work of S.~P.~M.~was partially supported by an NSERC Discovery Grant.
\end{acks}
\appendix

\section{Stage-value form}
\RK{} methods have an equivalent formulation in terms of stage values $\bfY^{(i)}$ that approximate the solution at time levels $t^n + c_i \Delta t$ rather than the derivatives.
If we start with an \RK{} tableau and apply this formulation to~\eqref{eq:odesystem}, we have
\begin{equation}
  \label{eq:value}
  \begin{split}
    \bfY^{(i)} & = \bfy^{(1,n)} + \Delta t \sum_{j=1}^s a_{ij} \bfZ^{(j)}, \\
    \bfZ^{(i)} & = \bfy^{(2,n)} + \Delta t \sum_{j=1}^s a_{ij} \bff(t^n+c_j \Delta t, \bfY^{(j)}, \bfZ^{(j)}),
  \end{split}
\end{equation}
where $\bfY^{(i)}$ and $\bfZ^{(i)}$ represent the function and its derivative, respectively.
If we substitute the first equation for $\bfY^{(i)}$ into the second, we can arrive at a stage equation for the $\bfZ^{(i)}$ alone:
\begin{equation}
  \bfZ^{(i)} = \bfy^{(2,n)} + \Delta t \sum_{j=1}^s a_{ij}
  \bff(t^n+c_j \Delta t, \bfy^{(1,n)} + \Delta t \sum_{k=1}^s a_{jk} \bfZ^{(k)} , \bfZ^{(j)}).
\end{equation}
This does not compress to a single tableau that admits Nystr\"om type methods like~\eqref{eq:nystrom_scheme}.

When there is no first derivative in the system, we can substitute the second equation in~\eqref{eq:value} into the first and arrive at a method for $\bfY^{(i)}$ alone:
\begin{equation}
  \begin{split}
    \bfY^{(i)} & = \bfy^{(1,n)} + \Delta t \sum_{j=1}^s a_{ij} \left( \bfy^{(2,n)} + \Delta t \sum_{k=1}^s a_{jk} \bff(t^n+c_k \Delta t, \bfY^{(k)}) \right) \\
    & = \bfy^{(1,n)} + \Delta t \sum_{j=1}^s a_{ij} \bfy^{(2,n)}
    + \left( \Delta t \right)^2 \sum_{j=1}^s a_{ij} \sum_{k=1}^s a_{jk} \bff(t^n+c_k \Delta t, \bfY^{(k)}) \\
    & = \bfy^{(1,n)} + c_i \left( \Delta t \right) \bfy^{(2,n)}
    + \left( \Delta t \right)^2 \sum_{k=1}^s \overline{a}_{ik} \bff(t^n+c_k \Delta t, \bfY^{(k)}) .
  \end{split}
\end{equation}

Since our first formulation does not support general \RKN{} schemes and our second is restricted to particular equations, we have not pursued these special cases in our implementation.

\bibliographystyle{ACM-Reference-Format}
\bibliography{references}


\begin{thebibliography}{46}


\ifx \showCODEN    \undefined \def \showCODEN     #1{\unskip}     \fi
\ifx \showDOI      \undefined \def \showDOI       #1{#1}\fi
\ifx \showISBNx    \undefined \def \showISBNx     #1{\unskip}     \fi
\ifx \showISBNxiii \undefined \def \showISBNxiii  #1{\unskip}     \fi
\ifx \showISSN     \undefined \def \showISSN      #1{\unskip}     \fi
\ifx \showLCCN     \undefined \def \showLCCN      #1{\unskip}     \fi
\ifx \shownote     \undefined \def \shownote      #1{#1}          \fi
\ifx \showarticletitle \undefined \def \showarticletitle #1{#1}   \fi
\ifx \showURL      \undefined \def \showURL       {\relax}        \fi
\providecommand\bibfield[2]{#2}
\providecommand\bibinfo[2]{#2}
\providecommand\natexlab[1]{#1}
\providecommand\showeprint[2][]{arXiv:#2}

\bibitem[\protect\citeauthoryear{Abu-Labdeh, MacLachlan, and
  Farrell}{Abu-Labdeh et~al\mbox{.}}{2023}]%
        {abu2022monolithic}
\bibfield{author}{\bibinfo{person}{Razan Abu-Labdeh}, \bibinfo{person}{Scott
  MacLachlan}, {and} \bibinfo{person}{Patrick~E. Farrell}.}
  \bibinfo{year}{2023}\natexlab{}.
\newblock \showarticletitle{Monolithic multigrid for implicit {R}unge-{K}utta
  discretizations of incompressible fluid flow}.
\newblock \bibinfo{journal}{\emph{J. Comput. Phys.}}  \bibinfo{volume}{478}
  (\bibinfo{year}{2023}), \bibinfo{pages}{111961}.
\newblock
\urldef\tempurl%
\url{https://doi.org/10.1016/j.jcp.2023.111961}
\showURL{%
\tempurl}


\bibitem[\protect\citeauthoryear{Adler, He, Hu, MacLachlan, and Ohm}{Adler
  et~al\mbox{.}}{2023}]%
        {JAdler_etal_2021a}
\bibfield{author}{\bibinfo{person}{James~H. Adler}, \bibinfo{person}{Yunhui
  He}, \bibinfo{person}{Xiaozhe Hu}, \bibinfo{person}{Scott MacLachlan}, {and}
  \bibinfo{person}{Peter Ohm}.} \bibinfo{year}{2023}\natexlab{}.
\newblock \showarticletitle{Monolithic multigrid for a reduced-quadrature
  discretization of poroelasticity}.
\newblock \bibinfo{journal}{\emph{SIAM J. Sci. Comput.}} \bibinfo{volume}{45},
  \bibinfo{number}{3} (\bibinfo{year}{2023}), \bibinfo{pages}{S54--S81}.
\newblock


\bibitem[\protect\citeauthoryear{Adler, Hu, MacLachlan, and Vijendiran}{Adler
  et~al\mbox{.}}{2025}]%
        {scott_biot}
\bibfield{author}{\bibinfo{person}{James~H. Adler}, \bibinfo{person}{Xiaozhe
  Hu}, \bibinfo{person}{Scott MacLachlan}, {and} \bibinfo{person}{Selvabavitha
  Vijendiran}.} \bibinfo{year}{2025}\natexlab{}.
\newblock \showarticletitle{Robust monolithic multigrid solvers for
  {R}unge--{K}utta discretizations of {B}iot elasticity}.
\newblock  (\bibinfo{year}{2025}).
\newblock
\newblock
\shownote{In preparation.}


\bibitem[\protect\citeauthoryear{Biot}{Biot}{1956}]%
        {10.1121/1.1908239}
\bibfield{author}{\bibinfo{person}{Maurice~A. Biot}.}
  \bibinfo{year}{1956}\natexlab{}.
\newblock \showarticletitle{Theory of Propagation of Elastic Waves in a
  Fluid‐Saturated Porous Solid. I. Low‐Frequency Range}.
\newblock \bibinfo{journal}{\emph{The Journal of the Acoustical Society of
  America}} \bibinfo{volume}{28}, \bibinfo{number}{2} (\bibinfo{date}{03}
  \bibinfo{year}{1956}), \bibinfo{pages}{168--178}.
\newblock
\urldef\tempurl%
\url{https://doi.org/10.1121/1.1908239}
\showDOI{\tempurl}


\bibitem[\protect\citeauthoryear{Biswas, Ketcheson, Seibold, and
  Shirokoff}{Biswas et~al\mbox{.}}{2023}]%
        {biswas2023design}
\bibfield{author}{\bibinfo{person}{Abhijit Biswas}, \bibinfo{person}{David~I.
  Ketcheson}, \bibinfo{person}{Benjamin Seibold}, {and} \bibinfo{person}{David
  Shirokoff}.} \bibinfo{year}{2023}\natexlab{}.
\newblock \showarticletitle{Design of {DIRK} schemes with high weak stage
  order}.
\newblock \bibinfo{journal}{\emph{Communications in Applied Mathematics and
  Computational Science}} \bibinfo{volume}{18}, \bibinfo{number}{1}
  (\bibinfo{year}{2023}), \bibinfo{pages}{1--28}.
\newblock


\bibitem[\protect\citeauthoryear{Both, Barnafi, Radu, Zunino, and
  Quarteroni}{Both et~al\mbox{.}}{2022}]%
        {BOTH2022114183}
\bibfield{author}{\bibinfo{person}{J.W. Both}, \bibinfo{person}{N.A. Barnafi},
  \bibinfo{person}{F.A. Radu}, \bibinfo{person}{P. Zunino}, {and}
  \bibinfo{person}{A. Quarteroni}.} \bibinfo{year}{2022}\natexlab{}.
\newblock \showarticletitle{Iterative splitting schemes for a soft material
  poromechanics model}.
\newblock \bibinfo{journal}{\emph{Computer Methods in Applied Mechanics and
  Engineering}}  \bibinfo{volume}{388} (\bibinfo{year}{2022}),
  \bibinfo{pages}{114183}.
\newblock
\showISSN{0045-7825}
\urldef\tempurl%
\url{https://doi.org/10.1016/j.cma.2021.114183}
\showDOI{\tempurl}


\bibitem[\protect\citeauthoryear{Brown, Knepley, May, McInnes, and Smith}{Brown
  et~al\mbox{.}}{2012}]%
        {brown2012composable}
\bibfield{author}{\bibinfo{person}{Jed Brown}, \bibinfo{person}{Matthew~G.
  Knepley}, \bibinfo{person}{David~A. May}, \bibinfo{person}{Lois~Curfman
  McInnes}, {and} \bibinfo{person}{Barry Smith}.}
  \bibinfo{year}{2012}\natexlab{}.
\newblock \showarticletitle{Composable linear solvers for multiphysics}. In
  \bibinfo{booktitle}{\emph{2012 11th International Symposium on Parallel and
  Distributed Computing}}. IEEE, \bibinfo{pages}{55--62}.
\newblock


\bibitem[\protect\citeauthoryear{Brubeck and Kirby}{Brubeck and Kirby}{2025}]%
        {brubeck2025fiat}
\bibfield{author}{\bibinfo{person}{Pablo~D. Brubeck} {and}
  \bibinfo{person}{Robert~C. Kirby}.} \bibinfo{year}{2025}\natexlab{}.
\newblock \showarticletitle{{FIAT}: enabling classical and modern
  macroelements}.
\newblock \bibinfo{journal}{\emph{arXiv preprint arXiv:2501.14599}}
  (\bibinfo{year}{2025}).
\newblock


\bibitem[\protect\citeauthoryear{Carpenter, Gottlieb, Abarbanel, and
  Don}{Carpenter et~al\mbox{.}}{1995}]%
        {carpenter1995theoretical}
\bibfield{author}{\bibinfo{person}{Mark~H Carpenter}, \bibinfo{person}{David
  Gottlieb}, \bibinfo{person}{Saul Abarbanel}, {and} \bibinfo{person}{Wai-Sun
  Don}.} \bibinfo{year}{1995}\natexlab{}.
\newblock \showarticletitle{The theoretical accuracy of {Runge--K}utta time
  discretizations for the initial boundary value problem: a study of the
  boundary error}.
\newblock \bibinfo{journal}{\emph{SIAM Journal on Scientific Computing}}
  \bibinfo{volume}{16}, \bibinfo{number}{6} (\bibinfo{year}{1995}),
  \bibinfo{pages}{1241--1252}.
\newblock


\bibitem[\protect\citeauthoryear{Chin-Joe-Kong, Mulder, and
  Van~Veldhuizen}{Chin-Joe-Kong et~al\mbox{.}}{1999}]%
        {chin1999higher}
\bibfield{author}{\bibinfo{person}{MJS Chin-Joe-Kong}, \bibinfo{person}{Wim~A
  Mulder}, {and} \bibinfo{person}{M Van~Veldhuizen}.}
  \bibinfo{year}{1999}\natexlab{}.
\newblock \showarticletitle{Higher-order triangular and tetrahedral finite
  elements with mass lumping for solving the wave equation}.
\newblock \bibinfo{journal}{\emph{Journal of Engineering Mathematics}}
  \bibinfo{volume}{35} (\bibinfo{year}{1999}), \bibinfo{pages}{405--426}.
\newblock


\bibitem[\protect\citeauthoryear{Clines, Howle, and Long}{Clines
  et~al\mbox{.}}{2022}]%
        {clines2022efficient}
\bibfield{author}{\bibinfo{person}{Michael~R Clines},
  \bibinfo{person}{Victoria~E Howle}, {and} \bibinfo{person}{Katharine~R
  Long}.} \bibinfo{year}{2022}\natexlab{}.
\newblock \showarticletitle{Efficient Order-Optimal Preconditioners for
  Implicit {Runge-Kutta} and {Runge-Kutta-Nystr\"{o}m} Methods Applicable to a
  Large Class of Parabolic and Hyperbolic {PDE}s}.
\newblock \bibinfo{journal}{\emph{arXiv preprint arXiv:2206.08991}}
  (\bibinfo{year}{2022}).
\newblock


\bibitem[\protect\citeauthoryear{Clough and Toucher}{Clough and
  Toucher}{1965}]%
        {clough1965finite}
\bibfield{author}{\bibinfo{person}{Ray~W. Clough} {and} \bibinfo{person}{J.~L.
  Toucher}.} \bibinfo{year}{1965}\natexlab{}.
\newblock \showarticletitle{Finite element stiffness matrices for analysis of
  plate bending}. In \bibinfo{booktitle}{\emph{Proc. of the First Conf. on
  Matrix Methods in Struct. Mech.}} \bibinfo{pages}{515--546}.
\newblock


\bibitem[\protect\citeauthoryear{Dahlquist}{Dahlquist}{1963}]%
        {dahlquist1963special}
\bibfield{author}{\bibinfo{person}{Germund~G. Dahlquist}.}
  \bibinfo{year}{1963}\natexlab{}.
\newblock \showarticletitle{A special stability problem for linear multistep
  methods}.
\newblock \bibinfo{journal}{\emph{BIT Numerical Mathematics}}
  \bibinfo{volume}{3}, \bibinfo{number}{1} (\bibinfo{year}{1963}),
  \bibinfo{pages}{27--43}.
\newblock


\bibitem[\protect\citeauthoryear{Dedner, Kloefkorn, and Nolte}{Dedner
  et~al\mbox{.}}{2020}]%
        {dedner2020python}
\bibfield{author}{\bibinfo{person}{Andreas Dedner}, \bibinfo{person}{Robert
  Kloefkorn}, {and} \bibinfo{person}{Martin Nolte}.}
  \bibinfo{year}{2020}\natexlab{}.
\newblock \showarticletitle{Python bindings for the DUNE-FEM module}.
\newblock \bibinfo{journal}{\emph{Zenodo. doi}}  \bibinfo{volume}{10}
  (\bibinfo{year}{2020}).
\newblock


\bibitem[\protect\citeauthoryear{Dupont}{Dupont}{1973}]%
        {dupont1973l2}
\bibfield{author}{\bibinfo{person}{Todd Dupont}.}
  \bibinfo{year}{1973}\natexlab{}.
\newblock \showarticletitle{$L^2$-estimates for Galerkin methods for second
  order hyperbolic equations}.
\newblock \bibinfo{journal}{\emph{SIAM J. Numer. Anal.}} \bibinfo{volume}{10},
  \bibinfo{number}{5} (\bibinfo{year}{1973}), \bibinfo{pages}{880--889}.
\newblock


\bibitem[\protect\citeauthoryear{Falgout and Yang}{Falgout and Yang}{2002}]%
        {falgout2002hypre}
\bibfield{author}{\bibinfo{person}{Robert~D Falgout} {and}
  \bibinfo{person}{Ulrike~Meier Yang}.} \bibinfo{year}{2002}\natexlab{}.
\newblock \showarticletitle{hypre: A library of high performance
  preconditioners}. In \bibinfo{booktitle}{\emph{Computational Science—ICCS
  2002: International Conference Amsterdam, The Netherlands, April 21--24, 2002
  Proceedings, Part III}}. Springer, \bibinfo{pages}{632--641}.
\newblock


\bibitem[\protect\citeauthoryear{Farrell, Kirby, and Marchena-Menendez}{Farrell
  et~al\mbox{.}}{2021a}]%
        {farrell2021irksome}
\bibfield{author}{\bibinfo{person}{Patrick~E. Farrell},
  \bibinfo{person}{Robert~C. Kirby}, {and} \bibinfo{person}{Jorge
  Marchena-Menendez}.} \bibinfo{year}{2021}\natexlab{a}.
\newblock \showarticletitle{{Irksome: Automating Runge--Kutta time-stepping for
  finite element methods}}.
\newblock \bibinfo{journal}{\emph{ACM Trans. Math. Software}}
  \bibinfo{volume}{47}, \bibinfo{number}{4} (\bibinfo{year}{2021}),
  \bibinfo{pages}{1--26}.
\newblock


\bibitem[\protect\citeauthoryear{Farrell, Knepley, Mitchell, and
  Wechsung}{Farrell et~al\mbox{.}}{2021b}]%
        {farrell2021pcpatch}
\bibfield{author}{\bibinfo{person}{Patrick~E. Farrell},
  \bibinfo{person}{Matthew~G. Knepley}, \bibinfo{person}{Lawrence Mitchell},
  {and} \bibinfo{person}{Florian Wechsung}.} \bibinfo{year}{2021}\natexlab{b}.
\newblock \showarticletitle{{PCPATCH}: software for the topological
  construction of multigrid relaxation methods}.
\newblock \bibinfo{journal}{\emph{ACM Transactions on Mathematical Software
  (TOMS)}} \bibinfo{volume}{47}, \bibinfo{number}{3} (\bibinfo{year}{2021}),
  \bibinfo{pages}{1--22}.
\newblock


\bibitem[\protect\citeauthoryear{Fischer}{Fischer}{1998}]%
        {fischer1998projection}
\bibfield{author}{\bibinfo{person}{Paul~F. Fischer}.}
  \bibinfo{year}{1998}\natexlab{}.
\newblock \showarticletitle{Projection techniques for iterative solution of
  {Ax= b} with successive right-hand sides}.
\newblock \bibinfo{journal}{\emph{Computer Methods in Applied Mechanics and
  Engineering}} \bibinfo{volume}{163}, \bibinfo{number}{1-4}
  (\bibinfo{year}{1998}), \bibinfo{pages}{193--204}.
\newblock


\bibitem[\protect\citeauthoryear{Frank, Schneid, and Ueberhuber}{Frank
  et~al\mbox{.}}{1985}]%
        {frank1985order}
\bibfield{author}{\bibinfo{person}{Reinhard Frank}, \bibinfo{person}{Josef
  Schneid}, {and} \bibinfo{person}{Christoph~W. Ueberhuber}.}
  \bibinfo{year}{1985}\natexlab{}.
\newblock \showarticletitle{Order results for implicit {R}unge--{K}utta methods
  applied to stiff systems}.
\newblock \bibinfo{journal}{\emph{SIAM J. Numer. Anal.}} \bibinfo{volume}{22},
  \bibinfo{number}{3} (\bibinfo{year}{1985}), \bibinfo{pages}{515--534}.
\newblock


\bibitem[\protect\citeauthoryear{Fu}{Fu}{2019}]%
        {fu2019-aa}
\bibfield{author}{\bibinfo{person}{Guosheng Fu}.}
  \bibinfo{year}{2019}\natexlab{}.
\newblock \showarticletitle{A high-order {HDG} method for the {B}iot's
  consolidation model}.
\newblock \bibinfo{journal}{\emph{Computers \& Mathematics with Applications}}
  \bibinfo{volume}{77} (\bibinfo{year}{2019}), \bibinfo{pages}{237--252}.
\newblock


\bibitem[\protect\citeauthoryear{Gee, Siefert, Hu, Tuminaro, and Sala}{Gee
  et~al\mbox{.}}{2006}]%
        {MWGee_etal_2006a}
\bibfield{author}{\bibinfo{person}{Michael~W. Gee}, \bibinfo{person}{Chris~M.
  Siefert}, \bibinfo{person}{Jonathan~J. Hu}, \bibinfo{person}{Ray~S.
  Tuminaro}, {and} \bibinfo{person}{Marzio~G. Sala}.}
  \bibinfo{year}{2006}\natexlab{}.
\newblock \bibinfo{booktitle}{\emph{{ML} 5.0 Smoothed Aggregation User's
  Guide}}.
\newblock \bibinfo{type}{{T}echnical {R}eport} SAND2006-2649.
  \bibinfo{institution}{Sandia National Laboratories}.
\newblock


\bibitem[\protect\citeauthoryear{Ham, Kelly, Mitchell, Cotter, Kirby, Sagiyama,
  Bouziani, Vorderwuelbecke, Gregory, Betteridge, Shapero, Nixon-Hill, Ward,
  Farrell, Brubeck, Marsden, Gibson, Homolya, Sun, McRae, Luporini, Gregory,
  Lange, Funke, Rathgeber, Bercea, and Markall}{Ham et~al\mbox{.}}{2023}]%
        {FiredrakeUserManual}
\bibfield{author}{\bibinfo{person}{David~A. Ham}, \bibinfo{person}{Paul H.~J.
  Kelly}, \bibinfo{person}{Lawrence Mitchell}, \bibinfo{person}{Colin~J.
  Cotter}, \bibinfo{person}{Robert~C. Kirby}, \bibinfo{person}{Koki Sagiyama},
  \bibinfo{person}{Nacime Bouziani}, \bibinfo{person}{Sophia Vorderwuelbecke},
  \bibinfo{person}{Thomas~J. Gregory}, \bibinfo{person}{Jack Betteridge},
  \bibinfo{person}{Daniel~R. Shapero}, \bibinfo{person}{Reuben~W. Nixon-Hill},
  \bibinfo{person}{Connor~J. Ward}, \bibinfo{person}{Patrick~E. Farrell},
  \bibinfo{person}{Pablo~D. Brubeck}, \bibinfo{person}{India Marsden},
  \bibinfo{person}{Thomas~H. Gibson}, \bibinfo{person}{Miklós Homolya},
  \bibinfo{person}{Tianjiao Sun}, \bibinfo{person}{Andrew T.~T. McRae},
  \bibinfo{person}{Fabio Luporini}, \bibinfo{person}{Alastair Gregory},
  \bibinfo{person}{Michael Lange}, \bibinfo{person}{Simon~W. Funke},
  \bibinfo{person}{Florian Rathgeber}, \bibinfo{person}{Gheorghe-Teodor
  Bercea}, {and} \bibinfo{person}{Graham~R. Markall}.}
  \bibinfo{year}{2023}\natexlab{}.
\newblock \bibinfo{booktitle}{\emph{Firedrake User Manual}
  (\bibinfo{edition}{first edition} ed.)}.
\newblock Imperial College London and University of Oxford and Baylor
  University and University of Washington.
\newblock
\urldef\tempurl%
\url{https://doi.org/10.25561/104839}
\showDOI{\tempurl}


\bibitem[\protect\citeauthoryear{Henson and Yang}{Henson and Yang}{2002}]%
        {VEHenson_UMYang_2002a}
\bibfield{author}{\bibinfo{person}{Van~Emden Henson} {and}
  \bibinfo{person}{Ulrike~Meyer Yang}.} \bibinfo{year}{2002}\natexlab{}.
\newblock \showarticletitle{Boomer{AMG}: a Parallel Algebraic Multigrid Solver
  and Preconditioner}.
\newblock \bibinfo{journal}{\emph{Applied Numerical Mathematics}}
  \bibinfo{volume}{41} (\bibinfo{year}{2002}), \bibinfo{pages}{155--177}.
\newblock


\bibitem[\protect\citeauthoryear{John, Linke, Merdon, Neilan, and Rebholz}{John
  et~al\mbox{.}}{2017}]%
        {john2017divergence}
\bibfield{author}{\bibinfo{person}{Volker John}, \bibinfo{person}{Alexander
  Linke}, \bibinfo{person}{Christian Merdon}, \bibinfo{person}{Michael Neilan},
  {and} \bibinfo{person}{Leo~G Rebholz}.} \bibinfo{year}{2017}\natexlab{}.
\newblock \showarticletitle{On the divergence constraint in mixed finite
  element methods for incompressible flows}.
\newblock \bibinfo{journal}{\emph{SIAM Rev.}} \bibinfo{volume}{59},
  \bibinfo{number}{3} (\bibinfo{year}{2017}), \bibinfo{pages}{492--544}.
\newblock


\bibitem[\protect\citeauthoryear{Kirby}{Kirby}{2024}]%
        {mmg}
\bibfield{author}{\bibinfo{person}{Robert~C Kirby}.}
  \bibinfo{year}{2024}\natexlab{}.
\newblock \showarticletitle{On the convergence of monolithic multigrid for
  implicit Runge-Kutta time stepping of finite element problems}.
\newblock \bibinfo{journal}{\emph{SIAM Journal on Scientific Computing}}
  (\bibinfo{year}{2024}).
\newblock
\newblock
\shownote{To appear.}


\bibitem[\protect\citeauthoryear{Kirby and MacLachlan}{Kirby and
  MacLachlan}{2025}]%
        {kirby2024extending}
\bibfield{author}{\bibinfo{person}{Robert~C Kirby} {and}
  \bibinfo{person}{Scott~P MacLachlan}.} \bibinfo{year}{2025}\natexlab{}.
\newblock \showarticletitle{Extending Irksome: improvements in automated
  Runge--Kutta time stepping for finite element methods}.
\newblock \bibinfo{journal}{\emph{ACM Transactions on Mathematical Software
  (TOMS)}} (\bibinfo{year}{2025}).
\newblock
\urldef\tempurl%
\url{https://doi.org/10.1145/3759245}
\showDOI{\tempurl}
\newblock
\shownote{To appear.}


\bibitem[\protect\citeauthoryear{Kirby and Mitchell}{Kirby and
  Mitchell}{2018}]%
        {kirby2018solver}
\bibfield{author}{\bibinfo{person}{Robert~C. Kirby} {and}
  \bibinfo{person}{Lawrence Mitchell}.} \bibinfo{year}{2018}\natexlab{}.
\newblock \showarticletitle{Solver composition across the {PDE}/linear algebra
  barrier}.
\newblock \bibinfo{journal}{\emph{SIAM Journal on Scientific Computing}}
  \bibinfo{volume}{40}, \bibinfo{number}{1} (\bibinfo{year}{2018}),
  \bibinfo{pages}{C76--C98}.
\newblock
\urldef\tempurl%
\url{https://doi.org/10.1137/17M1133208}
\showDOI{\tempurl}


\bibitem[\protect\citeauthoryear{Kraus, Lymbery, and Osthues}{Kraus
  et~al\mbox{.}}{2025}]%
        {kraus2024analysis}
\bibfield{author}{\bibinfo{person}{Johannes Kraus}, \bibinfo{person}{Maria
  Lymbery}, {and} \bibinfo{person}{Kevin Osthues}.}
  \bibinfo{year}{2025}\natexlab{}.
\newblock \showarticletitle{Time-continuous strongly conservative space-time
  finite element methods for the dynamic Biot model}.
\newblock  (\bibinfo{year}{2025}).
\newblock
\showeprint[arxiv]{2507.19955}~[math.NA]
\urldef\tempurl%
\url{https://arxiv.org/abs/2507.19955}
\showURL{%
\tempurl}


\bibitem[\protect\citeauthoryear{Logg, Mardal, and Wells}{Logg
  et~al\mbox{.}}{2012}]%
        {Logg:2012}
\bibfield{editor}{\bibinfo{person}{Anders Logg}, \bibinfo{person}{Kent-Andre
  Mardal}, {and} \bibinfo{person}{Garth~N. Wells}} (Eds.).
  \bibinfo{year}{2012}\natexlab{}.
\newblock \bibinfo{booktitle}{\emph{{Automated solution of differential
  equations by the finite element method: the FEniCS book}}}.
  Vol.~\bibinfo{volume}{84}.
\newblock \bibinfo{publisher}{Springer}.
\newblock
\urldef\tempurl%
\url{https://doi.org/10.1007/978-3-642-23099-8}
\showDOI{\tempurl}


\bibitem[\protect\citeauthoryear{Mitchell and M{\"u}ller}{Mitchell and
  M{\"u}ller}{2016}]%
        {mitchell2016high}
\bibfield{author}{\bibinfo{person}{Lawrence Mitchell} {and}
  \bibinfo{person}{Eike~Hermann M{\"u}ller}.} \bibinfo{year}{2016}\natexlab{}.
\newblock \showarticletitle{High level implementation of geometric multigrid
  solvers for finite element problems: {A}pplications in atmospheric
  modelling}.
\newblock \bibinfo{journal}{\emph{J. Comput. Phys.}}  \bibinfo{volume}{327}
  (\bibinfo{year}{2016}), \bibinfo{pages}{1--18}.
\newblock


\bibitem[\protect\citeauthoryear{Nilssen, Staff, and Mardal}{Nilssen
  et~al\mbox{.}}{2011}]%
        {nilssen2011order}
\bibfield{author}{\bibinfo{person}{Trygve~K Nilssen}, \bibinfo{person}{Gunnar~A
  Staff}, {and} \bibinfo{person}{Kent-Andre Mardal}.}
  \bibinfo{year}{2011}\natexlab{}.
\newblock \showarticletitle{Order optimal preconditioners for fully implicit
  {R}unge-{K}utta schemes applied to the bidomain equations}.
\newblock \bibinfo{journal}{\emph{Numerical Methods for Partial Differential
  Equations}} \bibinfo{volume}{27}, \bibinfo{number}{5} (\bibinfo{year}{2011}),
  \bibinfo{pages}{1290--1312}.
\newblock


\bibitem[\protect\citeauthoryear{Nystr\"{o}m}{Nystr\"{o}m}{1925}]%
        {nystrom1925}
\bibfield{author}{\bibinfo{person}{E.~J. Nystr\"{o}m}.}
  \bibinfo{year}{1925}\natexlab{}.
\newblock \showarticletitle{\"{U}ber die numerische Integration von
  Differentialgleichungen}.
\newblock \bibinfo{journal}{\emph{Acta Soc. Sci. Fenn.}} \bibinfo{volume}{50},
  \bibinfo{number}{13} (\bibinfo{year}{1925}), \bibinfo{pages}{1--54}.
\newblock


\bibitem[\protect\citeauthoryear{Pathria}{Pathria}{1997}]%
        {pathria1997correct}
\bibfield{author}{\bibinfo{person}{D Pathria}.}
  \bibinfo{year}{1997}\natexlab{}.
\newblock \showarticletitle{The correct formulation of intermediate boundary
  conditions for {Runge--K}utta time integration of initial boundary value
  problems}.
\newblock \bibinfo{journal}{\emph{SIAM Journal on Scientific Computing}}
  \bibinfo{volume}{18}, \bibinfo{number}{5} (\bibinfo{year}{1997}),
  \bibinfo{pages}{1255--1266}.
\newblock


\bibitem[\protect\citeauthoryear{Pazner and Persson}{Pazner and
  Persson}{2017}]%
        {pazner2017stage}
\bibfield{author}{\bibinfo{person}{Will Pazner} {and} \bibinfo{person}{Per-Olof
  Persson}.} \bibinfo{year}{2017}\natexlab{}.
\newblock \showarticletitle{Stage-parallel fully implicit {R}unge--{K}utta
  solvers for discontinuous {G}alerkin fluid simulations}.
\newblock \bibinfo{journal}{\emph{J. Comput. Phys.}}  \bibinfo{volume}{335}
  (\bibinfo{year}{2017}), \bibinfo{pages}{700--717}.
\newblock


\bibitem[\protect\citeauthoryear{Rana, Howle, Long, Meek, and Milestone}{Rana
  et~al\mbox{.}}{2021}]%
        {masud2021new}
\bibfield{author}{\bibinfo{person}{Md~Masud Rana}, \bibinfo{person}{Victoria~E
  Howle}, \bibinfo{person}{Katharine Long}, \bibinfo{person}{Ashley Meek},
  {and} \bibinfo{person}{William Milestone}.} \bibinfo{year}{2021}\natexlab{}.
\newblock \showarticletitle{A new block preconditioner for implicit
  {Runge--Kutta} methods for parabolic {PDE} problems}.
\newblock \bibinfo{journal}{\emph{SIAM Journal on Scientific Computing}}
  \bibinfo{volume}{43}, \bibinfo{number}{5} (\bibinfo{year}{2021}),
  \bibinfo{pages}{S475--S495}.
\newblock


\bibitem[\protect\citeauthoryear{Rani, Ghysels, Howle, Long, and Outrata}{Rani
  et~al\mbox{.}}{2025}]%
        {rani2025efficient}
\bibfield{author}{\bibinfo{person}{Aman Rani}, \bibinfo{person}{Pieter
  Ghysels}, \bibinfo{person}{Victoria Howle}, \bibinfo{person}{Katharine Long},
  {and} \bibinfo{person}{Michal Outrata}.} \bibinfo{year}{2025}\natexlab{}.
\newblock \showarticletitle{Efficient solution of fully implicit {Runge--K}utta
  methods for linear wave equations}.
\newblock \bibinfo{journal}{\emph{SIAM Journal on Scientific Computing}}
  (\bibinfo{year}{2025}), \bibinfo{pages}{S183--S206}.
\newblock


\bibitem[\protect\citeauthoryear{Rathgeber, Ham, Mitchell, Lange, Luporini,
  McRae, Bercea, Markall, and Kelly}{Rathgeber et~al\mbox{.}}{2016}]%
        {Rathgeber:2016}
\bibfield{author}{\bibinfo{person}{Florian Rathgeber},
  \bibinfo{person}{David~A. Ham}, \bibinfo{person}{Lawrence Mitchell},
  \bibinfo{person}{Michael Lange}, \bibinfo{person}{Fabio Luporini},
  \bibinfo{person}{Andrew T.~T. McRae}, \bibinfo{person}{Gheorghe-Teodor
  Bercea}, \bibinfo{person}{Graham~R. Markall}, {and} \bibinfo{person}{Paul
  H.~J. Kelly}.} \bibinfo{year}{2016}\natexlab{}.
\newblock \showarticletitle{{Firedrake: automating the finite element method by
  composing abstractions}}.
\newblock \bibinfo{journal}{\emph{ACM Trans. Math. Software}}
  \bibinfo{volume}{43}, \bibinfo{number}{3} (\bibinfo{year}{2016}),
  \bibinfo{pages}{24:1--24:27}.
\newblock
\urldef\tempurl%
\url{https://doi.org/10.1145/2998441}
\showDOI{\tempurl}
\showeprint[arxiv]{1501.01809}


\bibitem[\protect\citeauthoryear{Southworth, Krzysik, and Pazner}{Southworth
  et~al\mbox{.}}{2022a}]%
        {southworth2022fast1}
\bibfield{author}{\bibinfo{person}{Ben~S. Southworth}, \bibinfo{person}{Oliver
  Krzysik}, {and} \bibinfo{person}{Will Pazner}.}
  \bibinfo{year}{2022}\natexlab{a}.
\newblock \showarticletitle{Fast solution of fully implicit {R}unge-{K}utta and
  discontinuous {G}alerkin in time for numerical {PDE}s, {P}art {II}:
  {N}onlinearities and {DAE}s}.
\newblock \bibinfo{journal}{\emph{SIAM J. Sci. Comput.}} \bibinfo{volume}{44},
  \bibinfo{number}{2} (\bibinfo{year}{2022}), \bibinfo{pages}{A636--A663}.
\newblock
\showISSN{1064-8275}
\urldef\tempurl%
\url{https://doi.org/10.1137/21M1390438}
\showDOI{\tempurl}


\bibitem[\protect\citeauthoryear{Southworth, Krzysik, Pazner, and
  De~Sterck}{Southworth et~al\mbox{.}}{2022b}]%
        {southworth2022fast2}
\bibfield{author}{\bibinfo{person}{Ben~S. Southworth}, \bibinfo{person}{Oliver
  Krzysik}, \bibinfo{person}{Will Pazner}, {and} \bibinfo{person}{Hans
  De~Sterck}.} \bibinfo{year}{2022}\natexlab{b}.
\newblock \showarticletitle{Fast solution of fully implicit {R}unge-{K}utta and
  discontinuous {G}alerkin in time for numerical {PDE}s, {P}art {I}: {T}he
  linear setting}.
\newblock \bibinfo{journal}{\emph{SIAM J. Sci. Comput.}} \bibinfo{volume}{44},
  \bibinfo{number}{1} (\bibinfo{year}{2022}), \bibinfo{pages}{A416--A443}.
\newblock
\showISSN{1064-8275}
\urldef\tempurl%
\url{https://doi.org/10.1137/21M1389742}
\showDOI{\tempurl}


\bibitem[\protect\citeauthoryear{Staff, Mardal, and Nilssen}{Staff
  et~al\mbox{.}}{2006}]%
        {staff2006preconditioning}
\bibfield{author}{\bibinfo{person}{Gunnar~A. Staff},
  \bibinfo{person}{Kent-Andre Mardal}, {and} \bibinfo{person}{Trygve~K.
  Nilssen}.} \bibinfo{year}{2006}\natexlab{}.
\newblock \showarticletitle{Preconditioning of fully implicit Runge-Kutta
  schemes for parabolic {PDE}s}.
\newblock \bibinfo{journal}{\emph{Modeling, Identification, and Control}}
  \bibinfo{volume}{27}, \bibinfo{number}{1} (\bibinfo{year}{2006}),
  \bibinfo{pages}{109--123}.
\newblock


\bibitem[\protect\citeauthoryear{Van~Lent and Vandewalle}{Van~Lent and
  Vandewalle}{2005}]%
        {vanlent2005}
\bibfield{author}{\bibinfo{person}{Jan Van~Lent} {and} \bibinfo{person}{Stefan
  Vandewalle}.} \bibinfo{year}{2005}\natexlab{}.
\newblock \showarticletitle{Multigrid methods for implicit {Runge--Kutta} and
  boundary value method discretizations of parabolic {PDEs}}.
\newblock \bibinfo{journal}{\emph{SIAM Journal on Scientific Computing}}
  \bibinfo{volume}{27}, \bibinfo{number}{1} (\bibinfo{year}{2005}),
  \bibinfo{pages}{67--92}.
\newblock
\urldef\tempurl%
\url{https://doi.org/10.1137/030601144}
\showDOI{\tempurl}


\bibitem[\protect\citeauthoryear{Versbach, Linders, Kl{\"o}fkorn, and
  Birken}{Versbach et~al\mbox{.}}{2023}]%
        {versbach2023theoretical}
\bibfield{author}{\bibinfo{person}{Lea~Miko Versbach}, \bibinfo{person}{Viktor
  Linders}, \bibinfo{person}{Robert Kl{\"o}fkorn}, {and}
  \bibinfo{person}{Philipp Birken}.} \bibinfo{year}{2023}\natexlab{}.
\newblock \showarticletitle{Theoretical and practical aspects of space-time
  DG-SEM implementations}.
\newblock \bibinfo{journal}{\emph{The SMAI Journal of computational
  mathematics}}  \bibinfo{volume}{9} (\bibinfo{year}{2023}),
  \bibinfo{pages}{61--93}.
\newblock


\bibitem[\protect\citeauthoryear{Wanner and Hairer}{Wanner and Hairer}{1996}]%
        {wanner1996solving}
\bibfield{author}{\bibinfo{person}{Gerhard Wanner} {and} \bibinfo{person}{Ernst
  Hairer}.} \bibinfo{year}{1996}\natexlab{}.
\newblock \bibinfo{booktitle}{\emph{Solving Ordinary Differential Equations
  {II}}}.
\newblock \bibinfo{publisher}{Springer Berlin Heidelberg}.
\newblock


\bibitem[\protect\citeauthoryear{Wanner and N\o{}rsett}{Wanner and
  N\o{}rsett}{1993}]%
        {wanner1993solving}
\bibfield{author}{\bibinfo{person}{Gerhard Wanner} {and} \bibinfo{person}{Ernst
  N\o{}rsett, Syvert P.~and~Hairer}.} \bibinfo{year}{1993}\natexlab{}.
\newblock \bibinfo{booktitle}{\emph{Solving Ordinary Differential Equations
  {I}: Nonstiff Problems}}.
\newblock \bibinfo{publisher}{Springer Berlin Heidelberg}.
\newblock


\bibitem[\protect\citeauthoryear{zenodo/Zenodo-20250826.0}{zenodo/Zenodo-20250826.0}{2025}]%
        {zenodo/Zenodo-nirksome}
zenodo/Zenodo-20250826.0 \bibinfo{year}{2025}\natexlab{}.
\newblock \bibinfo{title}{{Software used in `Code for Automated RKN time
  stepping for finite element methods in Irksome'}}.
\newblock
\newblock
\urldef\tempurl%
\url{https://doi.org/10.5281/zenodo.16953980}
\showDOI{\tempurl}


\end{thebibliography}

\end{document}